\documentstyle[12pt,amssymb]{amsart}

\newcommand{\cR}{{\cal R}}

\newcommand{\cC}{{\cal C}}

\newcommand{\cL}{{\cal L}}
\newcommand{\cM}{{\cal M}}

\newcommand{\PGL}{\operatorname{PGL}}
\newcommand{\Vol}{\operatorname{Vol}}

\newcommand{\SO}{\operatorname{SO}}
\newcommand{\so}{\operatorname{so}}

\newcommand{\graph}{\operatorname{graph}}
        
        \newcommand{\bfmit}[1]{\hbox{\Bfmit {#1}}}

        \newfont{\Bfmit}{eufm10 scaled\magstep1}

\def\mapright#1{\smash{   \mathop{\longrightarrow }\limits^{#1}}}
\def\mapleft#1{\smash{   \mathop{\longleftarrow }\limits^{#1}}}

\def\proof{\smallskip\noindent{\it Proof. }}
\def\endproof{\hfill\qed}
\def\qed {\nobreak$\quad$\lower 1pt\vbox{
    \hrule
    \hbox to 8pt{\vrule height 8pt\hfil\vrule height 8pt}
      \hrule}\ifmmode\relax\else\par\medbreak \fi}

\newtheorem{thm}{Theorem}[section]

\newtheorem{lem}[thm]{Lemma}
\newtheorem{cor}[thm]{Corollary}

\newtheorem{prop}[thm]{Proposition}

\theoremstyle{definition}

\newtheorem{defn}[thm]{Definition}

\newtheorem{say}[thm]{}
\newtheorem{exmp}[thm]{Example}
\newtheorem{prob}[thm]{Problem}

\newtheorem{rem}[thm]{Remark}           

\theoremstyle{remark}

\begin{document}
\title[Moduli of stable polygons and $\overline{\cM}_{0,n}$]
{Moduli Spaces of stable Polygons and Symplectic Structures
on 
$\overline{\cM}_{0,n}$ 
}
\author[Yi {\sc Hu}]{Yi {\sc Hu}$^*$}

\thanks{
 To appear in {\bf Compositio Mathematica}} 
\maketitle

\begin{abstract}
In this paper, certain natural and elementary polygonal objects  
in Euclidean space, {\it the stable polygons}, are introduced, 
and the novel moduli spaces ${\bfmit M}_{{\bf r}, \varepsilon}$ of stable polygons
are constructed as complex analytic spaces. 
Quite unexpectedly,
these new moduli spaces are shown to be projective and isomorphic to 
the moduli space $\overline{\cM}_{0,n}$ of the Deligne-Mumford
stable curves of genus 0. Further, built into the structures of stable polygons
are some natural data
leading toward  to a family of (classes of) symplectic (K\"ahler) forms. 
To some degree, ${\bfmit M}_{{\bf r}, \varepsilon}$ may be considered as symplectic
counterparts of $\overline{\cM}_{0,n}$ and Kapranov's Chow quotient construction of
$\overline{\cM}_{0,n}$. All these together 
brings up a new tool to study 
the K\"ahler topology of $\overline{\cM}_{0,n}$.
\end{abstract}

{\footnotesize
\tableofcontents
}

\section{Introduction}

The moduli spaces of Riemann surfaces were introduced by Riemann in the nineteenth
century and have been since  great sources of interest and study. Most recently,
they have evolved as essential ingredients in symplectic topology and mirror symmetry.
In these theories, a guiding principle is that the geometry of an almost complex
manifold may be obtained by studying the space of all (pseudo-) holomorphic curves
on the manifold. One instance is the Gromov-Witten invariants. A key building block
in that theory is the Deligne-Mumford compactification $\overline{\cM}_{g,n}$ of the moduli space 
${\cM}_{g,n}$ of pointed Riemann surfaces of genus $g$. Symplectic topology and
mirror symmetry are, of course, not the only places where
$\overline{\cM}_{g,n}$ is important.
 These moduli spaces are also rich subjects of study
in Thurston's hyperbolic geometry and Teichm\"uller theory (and so on). On yet
another important aspect of geometry, the genus zero case 
$\overline{\cM}_{0,n}$ was recently effectively used by Kawamata to prove a
higher codimensional adjunction formula in algebraic geometry.

\smallskip
The basic questions we will address in this paper are:
{\sl how do we construct these compactified moduli spaces and
furthermore what geometric information can we derive
from the construction?} 

\smallskip

Our answers to these questions are very simple and rather unexpected.

\smallskip
We will only look at the special case of genus zero, namely, the
moduli space $\overline{\cM}_{0,n}$ of stable $n$-pointed Riemann spheres.
Set-theoretically, a generic point in  $\overline{\cM}_{0,n}$ is an 
equivalence class of  Riemann spheres with $n$ distinct ordered
points. All other points on the boundary of $\overline{\cM}_{0,n}$
are to be obtained by abiding to the following
principle: whenever some of marked points on a sphere come together,
we attach another sphere, called a bubble, at the coinciding position
and let the coinciding points get separated on the bubble sphere.
This process should go on until all  points get separated (see the picture below).

$$\overline{\cM}_{0,n} =$$

\vskip   2in
\begin{center}
\special{bmp:poly1.bmp x=4.8in y=2in}
Figure 1
\end{center}

In this paper, 
We will  (quite unexpectedly)
build $\overline{\cM}_{0,n}$ using totally differently materials. In short, 
we shall give {\it symplectic} constructions of $\overline{\cM}_{0,n}$
by using only elementary geometric  combinatorial
objects in the Euclidean space which
we will name as {\sl stable polygons}.
This yields some quite unexpected results which among other things
include a family of naturally and automatically built-in (classes of) symplectic K\"ahler forms
on $\overline{\cM}_{0,n}$. We believe that this will
eventually allow us to quantatively determine (in a forthcoming paper)
the K\"ahler cone and dually the Mori cone of effective curves
on $\overline{\cM}_{0,n}$ (cf. \S\S 8 and 9). 

\smallskip
Our approach is based upon a beautiful connection between
symplectic and algebraic geometry, discovered in late 70's and
early 80's. Briefly, let $(X, \omega)$
be a symplectic manifold that underlies a polarized projective
manifold $(X, \cL)$, where $\omega$ is a symplectic form on $X$
and $\cL$ is a positive line bundle such that $$c_1(\cL) = [\omega] \in H^2(X, {\Bbb Z})$$
Suppose a compact Lie group $K$ acts symplectically on $X$
with its complexification $G =K_{\Bbb C}$ acting holomorphically.
Then in a canonical way, we have a correspondence between quotients
in symplectic category and  quotients in projective algebraic category.
More precisely, assume that $K$ acts on $X$ in a Hamiltonian fashion,
that is, there is an equivariant moment map (depending on $\omega$)
$$\Phi_{\omega}: X \rightarrow {\bf k}^*$$
where ${\bf k}$ is the Lie algebra of $K$. Then the orbit space 
$\Phi_{\omega}^{-1}(0)/K$ carries a naturally induced symplectic structure
away from singularities. This is a symplectic quotient. On the other hand,
there is a dense open subset
$$X^{ss}(\cL) \subset X,$$
the set of all semistable points determined by $\cL$. By Mumford's
geometric invariant theory, $X^{ss}(\cL)$ has a quotient
$X^{ss}(\cL)/\!/G$ as a projective variety. The beautiful correspondence
mentioned above  is that
the symplectic quotient $\Phi_{\omega}^{-1}(0)/K$ is, in a canonical way,
homeomorphic to the projective quotient $X^{ss}(\cL)/\!/G$
(Mumford, Atiyah, Guillemin-Sternberg, Kempf-Ness, Ness, Kirwan, and others).

\smallskip
We now describe {\sl stable polygons} in ${\Bbb R}^3$,
the building blocks of our construction of $\overline{\cM}_{0,n}$.
A polygon in
${\Bbb R}^3$ is a collection of ordered vectors (edges)
$\{e_1, \ldots, e_n\}$ that add up to zero, $$e_1 + \ldots +e_n = 0.$$
Two polygons are equivalent if one can be obtained from the other
by a rigid Euclidean motion. Consider any $n$-gon with side lengths,
${\bf r} = (r_1, \ldots, r_n)$, then the deformation space of this $n$-gon
modulo equivalence is
$$M_{\bf r} =\{(e_1, \ldots, e_n) | e_1 + \ldots +e_n = 0, |e_i|=r_i, 1 \le i \le n \}/\sim.$$
This is a symplectic space with a symplectic form $\omega_{\bf r}$ away from
some (possibly) isolated singularities (for more about $M_{\bf r}$, 
see, e.g.,  \cite{KM96}. For the ring structure on $H^*(M_{\bf r})$, see \cite{HaKn97}).

\smallskip
To introduce stable polygons, we fix a term once and for all in this paper: 
edges are said to be {\it parallel} if they point the same direction; edges
pointing in opposite directions are considered {\it anti-parallel}, not parallel.

\smallskip

A generic polygon, that is, a polygon with no parallel edges, is stable in our sense. 
All other stable polygons are to be obtained
by abiding the following principle (much of which is 
 in the spirit of Fulton-MacPherson \cite{FM}): whenever a set of edges of a polygon
becomes parallel, we introduce  an (independent) generic
polygon, called a bubble,
 whose edges inherit the lengths of the abovementioned parallel edges
except the longest one. The longest length is the sum of the lengths
of the parallel edges minus a carefully chosen small positive
number $\epsilon$. (This $\epsilon$ has a precise quantity control and carries
significant symplecto-geometric meaning, consult Theorems 1.2 and 1.3 below.)
This process should go on until all sets of parallel edges
are properly addressed. 

\smallskip

In short, a stable polygon is a collection of labeled (but not ordered)
polygons that grows out of an ordinary polygon by introducing ``bubble''
polygons.  
The moduli space of all such stable polygons is denoted by
${\bfmit M}_{{\bf r}, \varepsilon}$. 
See the illustration below.
$${\bfmit M}_{{\bf r}, \varepsilon}=$$

\vskip   2.3in
\begin{center}
\special{bmp:poly2.bmp x=4.7in y=2.4in}

Figure 2
\end{center}

One should note that
the moduli space ${\bfmit M}_{{\bf r}, \varepsilon}$ of stable polygons
(in a very interesting manner)
 depends on two parameters: a fixed length vector ${\bf r}$ and
a collection $\varepsilon= \{\epsilon\}$ of   carefully chosen small positive
numbers. These built-in data in  stable polygons encode
significant symplecto-geometric information for $\overline{\cM}_{0,n}$.

\smallskip
Our main theorems are 

\begin{thm} {\rm (\S\S 5 and 7)}
The moduli space ${\bfmit M}_{{\bf r}, \varepsilon}$ has
a natural structure of a compact complex manifold and is
biholomorphic to the complex manifold that underlies the
projective variety $\overline{\cM}_{0,n}$.
\end{thm}


\begin{thm} {\rm (\S 6)}
By forgetting all the bubbles, we obtain a natural projection
$$\pi_{{\bf r}, \varepsilon}: {\bfmit M}_{{\bf r}, \varepsilon}
\rightarrow M_{\bf r}$$
that is holomorphic and bimeromorphic. It is the iterated blowup
of $M_{\bf r}$ along (the proper transforms of)
some explicitly described smooth subvarieties $Y_\alpha$ when
$M_{\bf r}$ is smooth. 
When $M_{\bf r}$ is singular, 
$\pi_{{\bf r}, \varepsilon}: {\bfmit M}_{{\bf r}, \varepsilon} \rightarrow M_{\bf r}$ is the composite of
a canonical resolution of singularities followed by explicit
iterated blowups.
\end{thm} 

In fact, more may be true. With the aid of the so-called symplectic $\varepsilon$-blowups
(e.g., see \cite{MS}), ${\bfmit M}_{{\bf r}, \varepsilon}$ (hence $\overline{\cM}_{0,n}$) comes equipped
with a symplectic form $\Omega_{{\bf r}, \varepsilon}$ whose cohomology class is uniquely defined
so that the map $$\pi_{{\bf r}, \varepsilon}: ({\bfmit M}_{{\bf r}, \varepsilon}, \Omega_{{\bf r}, \varepsilon})
\rightarrow (M_{\bf r}, \omega_{\bf r})$$
can be interpreted as a symplectic blowup. The ambiguity of the so-called symplectic blowups provides 
a natural and interesting  viewpoint to see
how natural it is the choices of ${\varepsilon}$ in the definition of stable polygons.
 This point calls for further investigation.

\smallskip

To close this introduction, 
some remarks are in order. Firstly,
$\overline{\cM}_{0,n}$ has a canonical K\"ahler class,
the Weil-Peterson class. We believe that this form should be realized
by some choice of $({\bf r}, \varepsilon)$. It has to come
 from regular $n$-gons. The choice of 
$\varepsilon$ is harder to guess (consult 4.6).
Secondly, by above, the collection of all  (legal\footnote{See 2.1 and 4.6 for the
legal ranges of ${\bf r}$ and $\varepsilon$.})
$({\bf r}, \varepsilon)$ should shed a light on the K\"ahler cone
of $\overline{\cM}_{0,n}$ (a conjecture is formulated in \S 9).
Finally, Thurston made remarks at author's seminar talk at UC, Davis that the same strategy
can be applied to unordered $n$-pointed Riemann spheres. Then the role
of polygons will then be replaced by collections of unordered vectors that
add up to zero. He also pointed out to the author that polyhedra may also be used
to study $\overline{\cM}_{0,n}$ (\cite{Thurston}).
Fried at author's talk at UC, Irvine,
also made similar suggestions about unordered $n$-pointed Riemann spheres.

\smallskip
Here is the structure of this paper. Sections 2 provides necessary
backgrounds on the symplectic geometry of the ordinary $n$-gons.
In Section 3, {\it favorable} moduli spaces of ordinary $n$-gons are
characterized in terms of their defining side length vectors (chambers). They are
useful throughout the paper. Section 4 introduces our principal
objects of study, stable polygons, mainly in set-theoretic terms. In Section 5, we show how
naturally stable polygons converge together to form a 
smooth compact complex analytic variety. 
Relation with
the moduli spaces of the ordinary polygons follows readily:
section 6 decomposes ${\bfmit M}_{{\bf r}, \varepsilon}$ as explicit
iterated blowups of $M_{\bf r}$.
Section 7 connects our new moduli spaces to  the
 moduli space $\overline{\cM}_{0,n}$ of
the Deligne-Mumford stable $n$-pointed curves of genus zero.
This leads to a new angle to look at the K\"ahler topology of
$\overline{\cM}_{0,n}$ which is discussed in the end.

\smallskip
There are some previous works on the moduli spaces of {\it ordinary} polygons
(\cite{W},  \cite{Kl}, \cite{KM96},  \cite{Kn}, among others) which may
fall into the general framework on relations between geometric invariant theory
 quotients and symplectic reductions 
(\cite{GIT}, \cite{Gui-Stern}, \cite{Kirwan}, \cite{Arms}, \cite{Hu91},
\cite{DH}, among others). The idea that $\overline{\cM}_{0,n}$ ought to have a counterpart
in terms of symplectic geometry was based on the observation (\cite{Kapranov95} and \cite{Hu91})
on a natural connection between GIT quotients and Chow (Hilbert) quotient.

\smallskip
{\bf Acknowledgments.} 
The author heartily thanks J\'anos Koll\'ar 
for providing the half of the financial support (Utah, 1995-96). 
Thanks are also due to Misha Kapovich for his 
collaboration during part of 1996 out of which much of this work grew
and John Millson for introducing him to the symplectic geometry of polygons.
He thanks Bill Thurston for his questions and comments and for mentioning
his related work \cite{Thurston}. He is grateful to Motohico Mulase
for inviting him to give a talk at UC, Davis, Greg Kuperberg and Albert Schwarz
for useful comments; Richard Wentworth and Mike Fried for inviting him to speak at UC, Irvine;
 Sean Keel for inviting a visit and talk at UT, Austin, and for his correction and comments;
 and the referee for his
many useful comments. The idea and part of this paper were started and written
when the author visited I.H.E.S, Bur-sur-Yvette (summer, 1995)
and Max-Planck-Institut, Bonn (summer, 1993, 96) whose hospitality and
financial supports are gratefully acknowledged. This work was made possible by
a Centennial Fellowship award by the American Mathematical Society.

\section{Geometry of ordinary polygons}

Many of the results on ordinary polygons 
that we shall collect below follows directly from  some well-known
general theory. See \cite{KM96}, for example, for a nice self-contained treatment.

\begin{say} A (spatial) polygon with $n$
sides in Euclidean space is determined by its labeled $n$
vertices $\{v_1, \cdots, v_n \}$ and these vertices are joined in cyclic order by
the directed edges $$\{e_1, \cdots, e_n \}$$ where $e_i$ starts from $v_i$ and ends at $v_{i+1}$
(here we set $v_{n+1} = v_1$). Since the $n$-gon is closed, 
$\{e_1, \cdots, e_n \}$,  regarded as vectors,
add up to zero, $$e_1 + \cdots + e_n = 0.$$
We consider the following equivalence relation among polygons.
Two $n$-gons are equivalent if one can be obtained from the other 
by the action of an orientation preserving
Euclidean isometry.  


Figure 3 illustrates three 12-gons. (The last two will be used later for other illustrations.)

\vskip   2in
\begin{center}
\special{bmp:poly3.bmp x=5in y=2in}

Figure 3
\end{center}
\end{say}

Throughout this paper, 
$M_{\bf r}$ will stand for the moduli space\footnote{The topology of this space is rather
easy to see. We will review some of its finer structures soon.} of $n$-gons with the 
       side length vector ${\bf r}=(r_1, \cdots, r_n)$. The normalized side length 
vector of ${\bf r}$,
        denoted by $\bar{{\bf r}}$, is defined to be $\frac{2} {L({\bf r})} {\bf r}$
        where $$L({\bf r}) = r_1 + \cdots + r_n$$ is the perimeter of the polygons
     of $M_{\bf r}$.
        As is well-known, $\bar{{\bf r}}$ lies in the hypersimplex
        $${\Bbb D}^n_2 := \{ (x_1, \ldots, x_n) \in {\Bbb R}^n : 0\le x_i \le 1, i =1, \ldots, n,
           \sum_i x_i = 2 \}.$$
To see the fine structures on $M_{\bf r}$,
we will appeal to the theory of symplectic reductions and K\"ahler geometric
invariant theory (GIT),  from which one will see that
        $M_{\bf r}$ and  $M_{\bar{{\bf r}}}$ are biholomorphic as complex analytic spaces
         and their symplectic (K\"ahler) forms (away from singularities)
        are proportional by the scalar  $\frac{2} {L}$.

\begin{say}
Consider the diagonal action of the group $\PGL (2, {\Bbb C})$ of projective transformations 
on $({\Bbb P}^1)^n$ and identify its maximal compact subgroup with $\SO (3)$. 
Identify ${\Bbb P}^1$ with the unit sphere $S^2$ in ${\Bbb R}^3$, $\SO (3)$-equivariantly. Let $\Vol (S^2)$ be the volume form on $S^2$.
Given any ${\bf r}$ with $r_1, \ldots, r_n > 0$,
$$\cL_{\bf r} = r_1 p_1^* \Vol (S^2) + \ldots + r_n p_n^*\Vol (S^2)$$
is a K\"ahler symplectic form on $({\Bbb P}^1)^n$ with respect to
 which the group $\SO (3)$
acts on $({\Bbb P}^1)^n$ in a Hamiltonian fashion. Here $p_i$ ($1 \le i \le n$) is the projection
from $({\Bbb P}^1)^n$ onto the $i$-th factor.
The moment map 
$$ \Phi_{\bf r}: ({\Bbb P}^1)^n \rightarrow \so (3)^* \cong {\Bbb R}^3 $$
determined by this Hamiltonian action is given as:
$$\Phi_{\bf r} (x_1, \ldots, x_n) = r_1 x_1  + \ldots + r_n x_n.$$

The correspondence between symplectic reductions
and K\"ahler (GIT, in the case that ${\bf r}$ is rational)
quotients asserts that
$$\Phi^{-1}_{\bf r} (0)/ \SO (3) \cong ({\Bbb P}^1)^n_{ss} ({\bf r}) /\!/ \PGL (2, {\Bbb C})$$
which is induced by the inclusion
$$\Phi^{-1}_{\bf r} (0) \subset  ({\Bbb P}^1)^n_{ss} ({\bf r}),$$
where $({\Bbb P}^1)^n_{ss} ({\bf r})$ denotes the set of $\PGL (2, {\Bbb C})$-semistable points
with respect to the K\"ahler form $\cL_{\bf r}$. We point out that
$\Phi^{-1}_{\bf r} (0) \ne \emptyset$ (or equivalently, 
$({\Bbb P}^1)^n_{ss} ({\bf r})\ne \emptyset$) if and only if $\bar{\bf r} \in {\Bbb D}^n_2$.
In this case, the quotient space 
$\Phi^{-1}_{\bf r} (0)/ \SO (3) \cong ({\Bbb P}^1)^n_{ss} ({\bf r}) /\!/ \PGL (2, {\Bbb C})$
has the expected dimension $n-3$ if and only if $\bar{\bf r} \in \hbox{int} {\Bbb D}^n_2$
(see \cite{Hu91}). For this reason, we make the assumption that
$\bar{\bf r} \in \hbox{int} {\Bbb D}^n_2$ throughout the rest of the paper.

For any polygon with edges $(e_1, \ldots, e_n)$
of the prescribed side lengths ${\bf r} = (r_1, \ldots, r_n)$
such that $\bar{\bf r} \in \hbox{int} {\Bbb D}^n_2$, set 
$$u_i = \frac{e_i}{r_i} \;\; \hbox{for} \;\; 1 \le i \le n.$$
Then $$(u_1, \cdots, u_n) \in (S^2)^n \cong ({\Bbb P}^1)^n$$
 and satisfies
$$r_1 u_i + \ldots + r_n u_n = e_1 + \ldots +  e_n = 0.$$
That is, $(u_1, \cdots, u_n) \in \Phi^{-1}_{\bf r} (0)$. Conversely, making the above arguments backward,
any point in $\Phi^{-1}_{\bf r} (0)$ determines a polygon with the prescribed side lengths ${\bf r} = (r_1, \ldots, r_n)$.
Thus we obtain the natural identification
$$M_{\bf r} \cong \Phi^{-1}_{\bf r} (0)/ \SO (3)$$
which, via the identification 
$$\Phi^{-1}_{\bf r} (0)/ \SO (3) \cong ({\Bbb P}^1)^n_{ss} ({\bf r}) /\!/ \PGL (2, {\Bbb C}),$$
 provides $M_{\bf r}$ a holomorphic structure and  a K\"ahler symplectic form $\omega_{\bf r}$
away from singularities. 
\end{say}

\begin{say}
\label{grasscons}
Alternatively, the moduli space $M_{\bf r}$ can also be constructed
as a symplectic quotient of the Grassmannian $G(2, {\Bbb C}^n)$ by the maximal
torus $({\Bbb C}^*)^n$. A simple way to see this is to look at the natural action
 of the group $({\Bbb C}^*)^n \times GL(2, {\Bbb C})$ on the space of full-rank matrices of size $2 \times n$ and then take quotients by stage. Using this approach to $M_{\bf r}$,
the results of \cite{DuiH} and \cite{GuiS89} (among others) can be applied directly.
\end{say}

\begin{say}
The hypersimplex is divided into chambers (maximal or otherwise,
making a polytopal chamber complex) by the walls defined as follows
$$W_J = \{ (x_1, \ldots, x_n) \in {\Bbb D}^n_2 :     \sum_{i \in J} x_i = 1 \}$$
where $J$ runs over all the proper subsets of $[n]=\{1, 2, \ldots, n\}$.
$J$ and its complement $J^c$ define the same wall $W_J=W_{J^c}$.
When $|J|=1$ or $n-1$, $W_J$ is a facet of ${\Bbb D}^n_2$ which is a {\it simplex}.
Other facets of ${\Bbb D}^n_2$ are of the form
$$F_i = \{ (x_1, \ldots, x_n) \in {\Bbb D}^n_2 : x_i = 0\}.$$
When $n > 4$, $F_i$ are again   hypersimpleces. 
A proper subset $J$ defines an interior wall if and only if $2 \le |J| \le n-2$ (i.e.,
$|J|>1$ and $|J^c|>1$).
Each interior wall $W_J$ divides ${\Bbb D}^n_2$ into two parts
$$\{ (x_1, \ldots, x_n) \in {\Bbb D}^n_2 : \sum_J x_j \le 1 \}$$ and
$$\{ (x_1, \ldots, x_n) \in {\Bbb D}^n_2 : \sum_{J^c} x_j \le 1 \}.$$
Two points $x$ and $y$ are in the same chamber (maximal or otherwise) if the following holds.
For all proper subset $J$,
$$\sum_J x_j \le 1 \Longleftrightarrow \sum_J y_j \le 1.$$
When only strict inequalities occur, we get characterizations of maximal chambers.
Points in the same chamber define homeomorphic (actually biholomorphic) moduli spaces of
polygons. But their naturally equipped structures of symplectic K\"ahler space
 are different, in general.
\end{say}

\begin{say}
We should also consider the positive cone over  ${\Bbb D}^n_2$ to take into account $n$-gons
with all possible side lengths
$$C({\Bbb D}^n_2) := \{ (x_1, \ldots, x_n) \in {\Bbb R}_+^n : \sum_i x_i \ne 0, 
           \frac{2(x_1, \ldots, x_n)}{\sum_i x_i} \in {\Bbb D}^n_2 \}$$
where ${\Bbb R}_+$ is the set of all nonnegative real numbers. Equivalently,
$$C({\Bbb D}^n_2) \cup \{0\} = \{ (x_1, \ldots, x_n) \in {\Bbb R}^n : 0 \le x_i \le \sum_{j \ne i} x_j,
i=1, \ldots, n\}.$$
For any $n$-gon $P \in M_{\bf r}$, the length vector  ${\bf r}$ belongs to $C({\Bbb D}^n_2)$, 
that is, $\bar{{\bf r}} \in {\Bbb D}^n_2$.
Walls and chambers in ${\Bbb D}^n_2$ induce obvious walls and chambers in
$C({\Bbb D}^n_2)$. Cone walls and cone chambers will generally be denoted by
$C(W_J)$ and $C(\bigtriangleup)$ (etc.) if $W_J$ is a wall of ${\Bbb D}^n_2$ and
$\bigtriangleup$ is chamber of ${\Bbb D}^n_2$.
All the above discussions extend to the cone $C({\Bbb D}^n_2)$ directly.
The dependence of the moduli spaces of
polygons on the cone chambers is the obvious one.
\end{say}

\begin{say}
\label{canisom}
Let ${\bf r}$ and ${\bf r}'$ be in the same chamber $C$ (maximal or otherwise).
$C$ is rational and thus contains a rational point ${\bf r}_C$. Let $M_C$ denote
the GIT quotient $({\Bbb P}^1)^n_{ss}({\bf r}_C)/\!/\PGL (2)$. This is a projective
variety and depends only on the chamber $C$.
Then we have the following isomorphisms
$${M}_{\bf r} \mapright{} {M}_C \mapleft{} {M}_{{\bf r}'}.$$
The two  isomorphisms are {\it canonical} because they are naturally induced
from the inclusions $\Phi_{\bf r}^{-1}(0) \subset ({\Bbb P}^1)^n_{ss}({\bf r}_C)
\supset \Phi_{{\bf r}'}^{-1}(0)$, 
which, of course, depend on no auxiliary
choices. Hence, we obtain an induced {\it canonical} isomorphism
$$f_{{\bf r}{\bf r}'}: M_{{\bf r}} \rightarrow M_{{\bf r}'}.$$
\end{say}

\begin{defn} A polygon is called a line gon if it lies entirely on a straight line.
\end{defn}

A line gon occurs if and only if there is a proper subset $J$ with
$$\sum_J r_i = \sum_{J^c} r_i .$$
When $n \ge 5$, they are  singular points of the moduli space and are isolated.

\begin{say}
\label{canmaps}
 Let $C$ be a chamber that lies on the boundary of another chamber $C'$.
By \cite{Hu91} (also \cite{DH} among others), there exists a canonical projective morphism
$$\beta_{C'C} : M_{C'} \rightarrow M_{C}$$
which, by \ref{canisom}, induces a canonical complex analytic map
$$\beta_{r'r} : M_{r'} \rightarrow M_{r}$$
where $r' \in C'$ and $r \in C$. When $C'$ is a maximal chamber,
the above maps are resolutions of singularities if $n > 4$.
\end{say}

\section{Favorable chambers}

\begin{say}
\label{favorchamber} For chambers adjacent to the boundary of the hypersimplex, 
the corresponding
moduli spaces take simple forms.

Let $\bigtriangleup_i$ be the (unique) maximal chamber in ${\Bbb D}^n_2$ that contains
the simplex facet $W_{\{i\}}$ ($1 \le i \le n$). Numerically, this chamber can be characterized
as follows. $ x \in \bigtriangleup_i$ if and only if 
$$ \sum_J x_j < 1$$ for all proper $J$ with $i \notin J$
and $|J| < n-1$.
Or equivalently (taking the complement $J^c$), 
$$ \sum_J x_j + x_i > 1$$
for all proper $J$ with $i \notin J$ if and only if
$$ x_j + x_i > 1$$ for all $j \ne i$.
\end{say}

\begin{defn} These chambers $\bigtriangleup_i$ ($i =1, \ldots, n$) and their
corresponding cones $C(\bigtriangleup_i)$ 
will play special r\^oles in this paper.
We will refer them as {\it favorable chambers}.
\end{defn}

\begin{prop} 
\label{propo:favorchamber}
(Theorem 2.1, \cite{Hu91} and Theorem 6.14, \cite{Hu96}) Let 
${\bf r}$ be an element in the interior of
the favorable chamber $C(\bigtriangleup_i)$.
Then  $M_{\bf r}$ is isomorphic to ${\Bbb P}^{n-3}$.
\end{prop}

\begin{prop} Assume that $n \ge 3$. Let ${\bf r} \in C(\bigtriangleup_i)$.
Then the $i$-th edge $e_i$ is the longest edge, that is, $r_i > r_j$ for all
$j \ne i$.
\end{prop}

\proof Dividing by $\frac{2}{L({\bf r})}$, we may assume that ${\bf r} \in \bigtriangleup_i$.
Assume otherwise that  $r_j \ge r_i$ for some $j$. Then
$$1 > \sum_{h \ne i} r_h = \sum_{h \ne i, j} r_h + r_j \ge \sum_{h \ne i, j} r_h + r_i > 1,$$
a contradiction.
\endproof

Recall from the introduction that for simplicity,
edges are said to be parallel if they point the
same direction. Edges pointing in opposite directions are considered to be
anti-parallel, not parallel.

\begin{defn}
\label{degeneration}
A polygon is said to degenerate at a set of
edges $e_I:=\{e_i\}_{i \in I}$ with $|I|>1$
if $\{e_i\}_{i \in I}$ are all parallel and no other
edges are parallel to them. An edge is said to be degenerate if it belongs
to a set of degenerate edges, $e_I$ with $|I|>1$.
\end{defn}

\begin{prop} 
\label{neverdeg} Let ${\bf r}$ be an element in the interior of
the favorable chamber $C(\bigtriangleup_i)$ and $P \in M_{\bf r}$.
Then the $i$-th edge of $P$  never degenerates along with others.
\end{prop}

\proof  Dividing by $\frac{2}{L({\bf r})}$, we may assume that ${\bf r} \in \bigtriangleup_i$.
Assume otherwise. Then
$$r_i + r_j \le \sum_{k \ne i, j} r_k$$ for some $j \ne i$ which implies that
$r_i + r_j \le 1$, a contradiction.
\endproof

\begin{rem}
\label{rem:p1n3}
Likewise, there is also a maximal chamber 
$\bigtriangledown_i$ in  ${\Bbb D}^n_2$ 
containing a (relatively) maximal chamber 
in the hypersimplex facet $F_{\{i\}}$ ($1 \le i \le n$) such that for any
${\bf r}$ in the interior of
the chamber $\bigtriangledown_i$, ${\cal M}_{\bf r}$ is isomorphic to $({\Bbb P}^1)^{n-3}$.
More precisely, a point $x$ of  ${\Bbb D}^n_2$ belongs to the chamber $\bigtriangledown_i$
if and only if 
$$ \sum_J x_j > 1$$ for all proper $J$ with $i \notin J$ if and only if
$$ x_j + x_k > 1$$ for all $j, k \ne i$.
Or equivalently (taking the complement $J^c$), 
$$ \sum_J x_j + x_i < 1$$
for all proper $J$ with $i \notin J$.
For more information about $\bigtriangledown_i$, consult Theorem 2.1 of \cite{Hu91}
and Theorem 6.10 of \cite{Hu96}.
We point out that these particular chambers will not play any special r\^oles in this paper
(although they may do if a different approach than the one in this paper is taken).
\end{rem}

\section{Stable degeneration of polygons}

To clarify the principal subject of our study, we shall give in this section a detailed set-theoretic
description of stable polygons, leaving out their finer structures to be treated in the
sequel.

\begin{defn} An $n$-gon is called generic if it does not have parallel edges.
Given a  side length vector ${\bf r}$,     we use $M_{\bf r}^0$ to denote the
       subspace of generic polygons.
\end{defn}

\begin{lem}
The moduli space $M_{\bf r}^0$ of all generic $n$-gons with the fixed side lengths
${\bf r}$ can be identified with the moduli space of $n$
distinct ordered points on the projective line.
\end{lem}

\proof This follows from the identification 
$M_{\bf r} \cong ({\Bbb P}^1)^n_{ss}({\bf r})/\!/\PGL (2)$.
\endproof

This is an open holomorphic space (in fact, quasi-projective)
and we are looking for  nice and geometrically meaningful compactifications
of it by adding   {\it limiting polygonal objects} such that the added objects fit together 
to form a divisor with normal crossings. 
The limiting polygonal objects that we choose to add  will be called stable polygons
which we will describe now.

\begin{say} 
The idea is the usual one as in Deligne-Mumford's construction of stable pointed curves.
We follow the spirit  of Fulton-MacPherson \cite{FM}.
The principal guiding principle is that whenever
a set of edges become parallel, we will resolve it by introducing a
''bubble'' polygon in a coherent way.
Here goes some  detailed descriptions of stable $n$-gons for $n \le 5$
to just gain some concrete feelings about general stable polygons.


All triangles are generic and thus stable.

A generic  quadrangle is stable. Given  a degenerate  quadrangle $P_0$ with
      two parallel edges $e_i$ and $e_j$,
      to remedy the ``coincidence'' problem,
    we introduce an independent arbitrary (stable) triangle $P_1$ with side 
     lengths $$(r_i, r_j, r_i+r_j -\epsilon_{i,j})$$
where $\epsilon_{i,j}$ is a fixed suitably\footnote{This will be made precise when
stable polygons are formally defined.} 
small positive number.
The collection $(P_0, P_1)$ of these two labeled (but not ordered) polygons is
a (reducible)  stable quadrangle.


 For pentagons, 
       one degeneration case is just like  quadrangles: 
       two edges $e_i$ and $e_j$ are parallel. In this case,
        we introduce an independent arbitrary (stable) triangle in exactly the same way 
        as we did in the quadrangle cases.

       When three edges, 
       $e_i$,  $e_j$, $e_k$, are parallel, to get a moduli-stable pentagon,
       we add an independent arbitrary moduli-stable quadrangle with prescribed side lengths
       $$(r_i, r_j, r_k,  r_i+r_j+r_k -\epsilon_{i,j,k})$$ where $\epsilon_{i,j, k}$
 is a suitably small  positive number.
The collection of these labeled (but not ordered) polygons   is a (reducible)  stable pentagon. 
\end{say}

\begin{say} In general,  given any polygon, if $k$ edges are parallel,
we add a generic $(k+1)$-gon whose first $k$ sides inherit the lengths
of the original edges but whose last side has a new length $r_ {i_1} + \cdots + r_{i_k} - 
\epsilon_{i_1, \ldots, i_k}$
with a choice of a small positive number $\epsilon_{i_1, \ldots, i_k}$. 
In any of the polygons obtained,
we allow more degeneration, and whenever edges become parallel,
we introduce (bubble) polygons as above. It is important to point out
 that in the further degeneration, the edges
with the new lengths\footnote{Such edges correspond to  double points of Deligne-Mumford stable curves,
as we shall see later.}
 will not degenerate (point the direction of
 any other edge). This will be automatic after choosing $\epsilon_{i_1, \ldots, i_k}$ carefully
(see Corollary \ref{rJinfc} below).
\end{say}

To formally define an arbitrary stable polygon, some preparations are in order.

\begin{lem}
\label{choosee}
Let ${\bf r}$ be any point in the interior of $C({\Bbb D}^n_2)$. Then
$$({\bf r}, \sum r_i - \epsilon)$$ is in the chamber of $C(\bigtriangleup_{n+1})$ of $C({\Bbb D}^{n+1}_2)$
if and only if $$\epsilon < 2 \min \{r_1, \ldots, r_n\}.$$
\end{lem}

\proof 
By \ref{favorchamber}, we need
$$\frac{2}{2 \sum r_i - \epsilon} \sum_J r_j < 1$$ for all proper subsets $J$ of $\{1, \ldots, n\}$, 
which is equivalent to
 $$\epsilon < 2 \sum_{J^c} r_i$$ for all proper subsets $J$ of $\{1, \ldots, n\}$,
which is, in turn,  equivalent to
$$\epsilon < 2 \min \{r_1, \ldots, r_n\}.$$
\endproof

\begin{say}
\label{chooseep}
Fix a vector ${\bf r}$ in the interior of the cone $C({\Bbb D}^n_2)$
over the hypersimplex ${\Bbb D}^n_2$.
Since ${\bf r}$ lies in the {\it interior} of the cone $C({\Bbb D}^n_2)$,
it is impossible to have $(n-1)$ many edges to point to the same direction.
It does not make sense to treat the case that ``there is one edge pointing to the
same direction''. Less trivial and a bit tricky is the case when
exactly two edges point  the same direction. According to the rule, we will have
to, in this case, introduce a triangle to make the polygon stable. However, 
because a triangle is rigid, it affects neither the moduli space
nor the symplectic structure (\S\S 5, 6 and 7). This last point will be important when
we consider the K\"ahler cone of $\overline{\cM}_{0,n}$.

So for any proper subset $J \subset \{1, \ldots, n\}$ with cardinality $1< |J| < n-1$
we  fix a (suitably small) positive number
$$0 < \epsilon_J < 2 \min_J  \{  r_j \}.$$
Set
$${\bf r}_{J, \epsilon_J} = (r_J, \sum_{J} r_j - \epsilon_J).$$
We frequently write ${\bf r}_{J, \epsilon_J}$ as ${\bf r}_{J}$ when no confusion may occur.
\end{say}

\begin{cor}
\label{rJinfc} The vector ${\bf r}_{J, \epsilon_J}$ lies in a favorable chamber.
Consequently, $M_{{\bf r}_{J}}$ is isomorphic to some projective space and
hence is independent of the choice of $\epsilon_J$.
Furthermore, no polygon in $M_{{\bf r}_{J}}$ ever degenerates at the last
edge (i.e., the longest edge).
\end{cor}

\begin{rem}
\label{aboute}
Note that for any $0< \epsilon < 2 \min \{r_1, \ldots, r_n\}$,  we can set $\epsilon_J = \epsilon$
for all proper subsets $J \subset [n]$. This will satisfy part of our purposes (cf. Corollary
\ref{rJinfc}) and may save some
notational mess. In particular, there exists a {\it canonical} choice
$$\epsilon_J = \min\{r_1, \ldots, r_n\}$$
for all proper subsets $J \subset [n]$. To keep generality, however, we will work with
an arbitrary choice of  $\epsilon_J$ in the legal range.
\end{rem}

\begin{rem}
We will use $C_J$ to denote the favorable chamber that ${\bf r}_J$ belongs to. In particular,
      $M_{C_J}$ stands for the common projective model for all $M_{\bf x}$ with
       ${\bf x} \in C_J$. Note that this chamber does not depend on the choice of $\epsilon_J$.
\end{rem}

\begin{defn} 
\label{bubblerelation}
Let ${\bf r}$ be a point in the interior of $C({\Bbb D}^n_2)$. Given any proper
subset $J$ of $\{1, \ldots, n\}$ with $|J| \ge 2$, a pair
$$(P, P') \in M_{\bf r} \times M_{{\bf r}_J}$$
is said to be a bubble pair if $P$ degenerates at the edges $e_J$ (i.e.,
the edges $\{e_j\}_J$ of $P$ point to the same direction
and no other edges point to the direction of these edges). In this
case, $P'$ is called a bubble of $P$. 
\end{defn}

\vskip   2.2in
\begin{center}
\special{bmp:poly4.bmp x=5in y=2.2in}
Figure 4
\end{center}

By Corollary \ref{rJinfc}, $P'$ never degenerates
at the last edge (i.e., the longest edge whose length is $\sum_J r_j - \epsilon_J$).
We point out that $(P, P')\in M_{\bf r} \times M_{{\bf r}_J}$
being a bubble pair implies that 
$$\sum_J r_j \le \sum_{J^c} r_j.$$  

\begin{defn}
For this reason, we shall call
such a proper subset {\it relevant}. Thus only {\it relevant} proper subset $J$
will occur as the index sets for bubble pairs. 
\end{defn}

Thus the relevancy is a relative concept.
 This point will be useful in the definition of
stable polygons. For simplicity, the set of all  relevant subsets $J \subset \{1, \ldots, n\}$
with respect to ${\bf r}$  is denoted by $\cR ({\bf r})$.

\begin{say}
Now let ${\bf r}$ be any point in the interior of $C({\Bbb D}^n_2)$ and
$${\varepsilon} = \{\epsilon_J | J \in \cR ({\bf r}) \}$$ with $\epsilon_J$ chosen as 
in \ref{chooseep}.
Fix them once and for all. (We point out again that by Remark \ref{aboute},
we may choose all $\epsilon_J$ equal to a fixed number $\epsilon$. And there is even an
canonical choice of such $\epsilon$, namely, $\min\{r_1, \ldots, r_n\}$.)

For a notational clarification, we make a remark that in this paper, a plain Greek letter
$\epsilon$ is always to mean a suitably small positive number, while ${\varepsilon}$
is to stand for a collection of such $\epsilon$.
\end{say}

\begin{defn} 
\label{definestablegons}
A stable $n$-gon with respect to the side length vector  ${\bf r}$ and 
${\varepsilon} = \{\epsilon_J | J \in \cR ({\bf r}) \}$
is a collection of labeled (but {\it not ordered}) polygons 
$${\bf P}:=(P_0, P_1, \cdots, P_m) \in {M}_{{\bf r}} \times {M}_{{\bf r}_{J_1}} \times \cdots 
   \times {M}_{{\bf r}_{J_m}}$$
satisfying the following properties:
\begin{enumerate}
\item whenever $J_t \subset J_s$, then $P_t$ is a bubble of $P_s$. 
   \item if $P_h$ does not have a bubble, then it is generic (i.e.,$P_h \in M^0_{{\bf r}_{J_h}}$). 
\end{enumerate}
\end{defn}

In particular, in Definition \ref{definestablegons} (1),
we must have that   $J_t \subset J_s$ is relevant with respect to
${\bf r}_{J_s}$.

One thing worths a word now. When $J \in \cR ({\bf r})$ and $|J|=2$, 
$M_{r_J}$ consists of a single triangle. Although it is convenient to include
it to define a stable polygon, in effect, it will not do anything to the
moduli space. In particular, when we later formulate the K\"ahler cone of
$\overline{\cal M}_{0,n}$, these ''small'' $J$'s have to be dropped out from computations.
For later use, we set $\cR_{>2}({\bf r}) = \{J \in \cR ({\bf r}) | |J|>2\}$.

\begin{rem} {\it
Using a permutation among edges, we can always arrange the edges pointing to the same
direction to be adjacent. The inverse of the permutation
will allow us to get back what we start with. }
This observation will  simplify our  exposition at points.
\end{rem}

Figure 5 illustrates two polygons that are related by the permutation
$$(1\; 2\; 3\;4\;5\;6\;7\;8\;9\;10\;11\;12) \to (1\;7\;2\;10\;3\;11\;4\;12\;5\;8\;6\;9).$$
The edges of the first polygon are numbered starting
from the left-most edge with the direction $\rightarrow$. 
 The edges of the second polygon are numbered in the same way.

\begin{center}
\vskip   2.1in
\special{bmp:poly5.bmp x=5.5in y=2.1in}

 Figure 5
\end{center}

\begin{defn}
We use ${\bfmit M}_{{\bf r}, {\varepsilon}}$ to denote the set of all
stable polygons. For simplicity, one may abbreviate 
${\bfmit M}_{{\bf r},  {\varepsilon}}$ 
as  ${\bfmit M}_{\bf r}$.
\end{defn}

It is rather easy to see that 
the set ${\bfmit M}_{{\bf r},{\varepsilon}}$
 of all stable polygons with the prescribed side vector 
${\bf r}$ and a choice of ${\varepsilon}$
carries a natural compact Hausdorff topology. This topology will be the underlying
topology of the complex structure on  ${\bfmit M}_{{\bf r}, {\varepsilon}}$
 which we will construct
in the next section.

\section{Moduli spaces ${\bfmit M}_{{\bf r}, \varepsilon}$ of stable polygons}

Throughout, unless stated otherwise,
${\bf r}$ denotes a fixed length vector away from the boundary of $C({\Bbb D}^n_2)$
and $C$ denotes the chamber (not necessarily maximal) that contains ${\bf r}$.

\begin{say}
Consider the space of $n$-polygons with one free side:
$$Z_{[n]} := \left\{ 
 \begin{array}{ccccc}
\hbox{$n$-gons with the non-zero fixed side length} \; r_1, \ldots, r_{n-1} \\
\hbox{ but the last side is free}
\end{array}
\right\}.$$
For any polygon $P$ in $Z_{[n]}$, define $l(P)$ to be the length of the free side
(i.e., the last side). 
Then this length function 
$$l: Z_{[n]} \rightarrow {\Bbb R}$$
assumes the maximal value $r_1 + \ldots + r_{n-1}$.
For convenience, we may allow the length of the free side to be zero. Then $l$ also assumes the
minimal value 0.

Every level set of $l$ provides a moduli space of polygons with the obvious prescribed
side lengths. At the extremal, $f^{-1}(0)$ is the moduli space of $(n-1)$-gons with
the prescribed  side lengths $( r_1, \ldots, r_{n-1})$; while 
$f^{-1}(r_1 + \ldots + r_{n-1})$ consists of a single line polygon.

We are interested in $M_{\bf r}= l^{-1}(r_n)$ where the length of the last side $r_n$
is close to the maximum $r_1 + \ldots + r_{n-1}$.
\end{say}

By Lemma \ref{choosee},
 if 
$$r_1 + \ldots + r_{n-1} - 2 \min\{r_1, \ldots, r_{n-1} \} < a < b < r_1 + \ldots + r_{n-1},$$
then $l^{-1}(a)$ and $l^{-1}(b)$ are both isomorphic to the same projective space. We may use
$f_{ab}$ to denote the {\it canonical} isomorphism from $l^{-1}(a)$ to $l^{-1}(b)$:
$$f_{ab}: l^{-1}(a) \rightarrow l^{-1}(b).$$
When $b=r_1 + \ldots + r_{n-1}$, $l^{-1}(b)$ is a single point. In this case,
we understand $f_{ab}$ as the total collapsing map $l^{-1}(a) \rightarrow l^{-1}(b)$.

\begin{say}
\label{say:incidentrelation}
Let $P_0 \in M_{\bf r}$ be an $n$-gon
and $J =\{j_1 < \ldots  < j_s \}$ a subset of $\{1, \ldots, n\}$ 
such that $P_0$ degenerates at $e_J$
(see Definition \ref{degeneration}). For any $P \in M_{\bf r}$, {\it via a permutation
if necessary, 
we can assume that the edges $e_{j_1}, \ldots, e_{j_s}$ are adjacent in that order}.
Now let 
$$d_J = -(e_{j_1} + \ldots + e_{j_s})$$
be the vector  starting from the end of $e_{j_s}$
and ending at the initial of $e_{j_1}$. 
 $d_J$ will be referred as the $J$-diagonal and can be zero in general.
Then $e_j (j \in J)$ and $d_J$ form a polygon
$Q_J$, while via a similar arrangement,
$-d_J$ and $e_j (j \in J^c)$ form a polygon $Q_J^c$. Informally, we may say that 
the $J$-diagonal $d_J$ divides the (permuted) $n$-gon $P$ 
into the union of two sub-polygons $Q_J$ and $Q_J^c$: $P := Q_J \cup Q_J^c$.
For our purpose, we assume that $l(Q_J)$ (=the length $|d_J|$ of the $J$-diagonal =
the length of the last side of $Q_J$)  is close to
the maximal value of $l$ in the sense that 
$$r_{j_1} + \ldots + r_{j_s} - 2 \min\{r_{j_1}, 
\ldots, r_{j_s} \} < l(Q_J) \le r_{j_1} + \ldots + r_{j_s}.$$
Then $(r_{j_1}, \ldots, r_{j_s}, l(Q_J))$ is in the favorable chamber $C_J$ 
unless $l(Q_J) = r_{j_1} + \ldots + r_{j_s}$ in which case
the normalization of $(r_{j_1}, \ldots, r_{j_s}, l(Q_J))$ 
belongs to a boundary simplex $C_J'$ in ${\Bbb D}^n_2$.
\end{say}

\begin{rem} Applying the inverse of the permutation (if necessary), we can remove
the assumption that edges $e_J$ are adjacent. For simplicity, rather than declaring
this at every place where we use it, we will simply say that ``the permutation scheme
is applied'' or ``apply the  permutation scheme''.
\end{rem}

\begin{defn} Let  $P_0 \in M_{\bf r}$ be an $n$-gon 
and $J =\{j_1 < \ldots  < j_s \}$ a subset of $\{1, \ldots, n\}$
such that $P_0$ degenerates at $e_J$. 
Let $P \in M_{\bf r}$ and $Q \in M_{{{\bf r}_J}}$. Let also
$Q_J$ and $Q_J^c$ be as in \ref{say:incidentrelation}. Here the permutation scheme
is applied when necessary.
We say that $P$ and $Q$ are incident,
denoted by $$\iota(P, Q),$$ 
if $$r_{j_1} + \ldots + r_{j_s} - 2 \min\{r_{j_1}, \ldots, r_{j_s} \} < 
l(Q_J)=|d_J| \le r_{j_1} + \ldots + r_{j_s}$$
and $Q$ is mapped to $Q_J$ by the morphism $f_{l(Q)l(Q_J)}$. 
The map $f_{l(Q)l(Q_J)}$ is an isomorphism unless $l(Q_J)$ is maximal (i.e., $Q_J$ is a line gon)
in which case $f_{l(Q)l(Q_J)}$ is a total collapsing map to a point variety.
\end{defn}

Figure 6 depicts an example where $Q$ and $Q_J$ are assumed to be
isomorphic by the morphism $f_{l(Q)l(Q_J)}$. 

\begin{center}
\vskip   2.5in
\special{bmp:poly6.bmp x=5.8in y=2.5in}

 Figure 6
\end{center}

Note that we automatically have for any $Q \in M_{{{\bf r}_J}}$ that 
$$r_{j_1} + \ldots + r_{j_s} - 2 \min\{r_{j_1}, \ldots, r_{j_s} \} < 
l(Q) < r_{j_1} + \ldots + r_{j_s}$$ because $l(Q) = \sum_J r_j - \epsilon_J$, by 
definition. Thus, 
the above definition makes sense because of the discussion in \ref{say:incidentrelation}.

\begin{lem}
\label{iisanalytic}
The incidence relation is complex analytic. That is, the subset of
$M_{\bf r} \times M_{{{\bf r}_J}}$ defined by incidence relation is
complex analytic.
\end{lem}

\proof
It basically follows from the fact that all maps $f_{l(Q')l(Q)}$
in the definition of incidence relation
are complex analytic. 
\endproof

\begin{say} Let $P_0 \in M_{\bf r}$ be a polygon that degenerates at the edges $e_J$.
Here, if necessary, we may apply the permutation scheme.
Set $U_{\bf r}(P_0)_J \subset M_{\bf r}$ to be the  subset of $M_{\bf r}$ consisting of
polygons $P$ such that 
$$r_{j_1} + \ldots + r_{j_s} - 2 \min\{r_{j_1}, \ldots, r_{j_s} \} < |d_J|
\le  r_{j_1} + \ldots + r_{j_s}. $$ This is an open neighborhood of $P_0$ in $M_{\bf r}$.
\end{say}

\begin{lem}
\label{nolinegons} We have
 \begin{enumerate}
\item If $P_0$ is not a line gon, $U_{\bf r}(P_0)_J$ does not contain a line gon. 
      In particular, $U_{\bf r}(P_0)_J$ is smooth.
\item If $P_0$ is a line gon, $P_0$ is the only line gon in $U_{\bf r}(P_0)_J$.
     In particular, $U_{\bf r}(P_0)_J$ is smooth away from $P_0$.
\end{enumerate}
\end{lem}

\proof We apply the permutation scheme and assume that the edges $e_J$ are adjacent.
Let $P_0'$ be a line gon in  $U_{\bf r}(P_0)_J$ that degenerates at edges $e_I$ for some $I \ne J$.
(Here, whether the edges $e_I$ are adjacent is irrelevant to our proof.)
In particular, we have $\sum_I r_i = \sum_{I^c} r_i$. Let $d_J$ be the $J$-diagonal for $P_0'$.
Without loss of generality, assume that $\sum_{I \cap J} r_j \ge \sum_{I^c \cap J} r_j $.
Then $$|d_J| = | \sum_{I \cap J} r_j - \sum_{I^c \cap J} r_j | = 
\sum_J r_j - 2 \sum_{I^c \cap J} r_j \le \sum_J r_j - 2 \min_J \{r_j\}.$$
This contradicts that  $P_0' \in U_{\bf r}(P_0)_J$.
\endproof

We will need a general lemma

\begin{lem} {\rm (Removable Singularity Theorem).}
\label{extendingmap}
Let $f: X \rightarrow Y$ be a continuous map between two holomorphic complex varieties.
Assume that  the restriction $f_0$ of $f$ to a dense open subset of $X$ is holomorphic. Then
$f$ is itself holomorphic.
\end{lem}

\proof
It follows from the consideration of the following diagram
$$f: X \rightarrow {\rm Graph}(f) = \overline{{\rm Graph}(f_0)}  \rightarrow X \times Y 
 \rightarrow Y.$$
\endproof

When $P_0$ is a line gon that degenerate at $e_J$ and $e_{J^c}$, then $e_J$ and $e_{J^c}$
are the only sets of degenerating edges for $P_0$. Set
$U^0_{\bf r}(P_0)= U_{\bf r}(P_0)_J \cap U_{\bf r}(P_0)_{J^c}$.
This is again an open (singular) neighborhood of $P_0$ in $M_{\bf r}$.
Define a correspondence set
$$\widehat{U}_{\bf r}^1({P_0}) \subset  U^0_{\bf r}(P_0) \times M_{{\bf r}_J} \times
M_{{\bf r}_{J^c}}$$ as follows
$$\widehat{U}_{\bf r}^1({P_0})=\{(P,P_1,P_2) \in M_{\bf r} \times M_{{\bf r}_J} \times
M_{{\bf r}_{J^c}} : \iota(P,P_1), \iota(P,P_2). \}$$
$\widehat{U}_{\bf r}^1({P_0})$ projects onto
$U^0_{\bf r}(P_0)$ by forgetting the last two factors.

\begin{lem}
\label{canonicalresolution}
Let $P_0$ be a line gon. Then the projection
$$\widehat{U}_{\bf r}^1({P_0}) \rightarrow U^0_{\bf r}({P_0})$$ is a (canonical) resolution
of singularities\footnote{This lemma is a special case of a more 
general result from the Geometric Invariant Theory or the theory of
symplectic reductions. Our proof is independent and uses only polygons.} of $U^0_{\bf r}(P_0)$.
\end{lem}

\proof We will apply the permutation scheme in this proof.
$\widehat{U}_{\bf r}^1({P_0})$ is the subset of
$U^0_{\bf r}({P_0}) \times M_{{\bf r}_J} \times M_{{\bf r}_{J^c}}$ 
 defined by the incidence relation.
It projects surjectively onto
the subset $U_{\bf r}(J^c)$ of $U^0_{\bf r}({P_0}) \times M_{{\bf r}_{J^c}}$ by forgetting
the middle factor. $U_{\bf r}(J^c)$ is nothing but
$$U_{\bf r}(J^c)=\{(P,P_1) \in U^0_{\bf r}({P_0}) \times M_{{\bf r}_{J^c}}: \iota(P,P_1)\}.
$$
We will first argue that the subset $U_{\bf r}(J^c)$ is smooth by identifying it
with a smooth open subset of $M_{{\bf r}'}$ where ${\bf r}'$ is a 
perturbation of ${\bf r}$.
 For this purpose, assume that 
$${\bf r}_J =(r_{j_1}, \ldots,  r_{j_s}, \sum_J r_j - \epsilon_J),\;
\hbox{and} \; {\bf r}_{J^c}=(r_{J^c}, \sum_{J^c} r_j - \epsilon_{J^c}).$$ Note that
$$\sum_J r_j = \sum_{J^c} r_j.$$
Let ${\bf r}_J'=(r_{j_1}', \ldots,  r_{j_s}', \sum_J r_j - \epsilon) \in C_J$
be a small generic  perturbation of 
${\bf r}_J= (r_{j_1}, \ldots,  r_{j_s}, \sum_J r_j - \epsilon_J)$ in $C_J$
such that 
$$r_{j_1}' + \ldots +  r_{j_s}' =  \sum_{J^c} r_j - \epsilon_{J^c}.$$
Let ${\bf r}'$ to be obtained from ${\bf r}$ by replacing $r_{j_1}, \ldots,  r_{j_s}$ by
$r_{j_1}', \ldots,  r_{j_s}'$ correspondingly.
Since ${\bf r}'$ is in more general position, and  can be taken
to be close to ${\bf r}$, we get a canonical surjection
$\beta: M_{{\bf r}'} \rightarrow  M_{\bf r}$ by \ref{canmaps}.

Let $U_{\bf r'} = \beta^{-1}(U^0_{\bf r}(P_0))$. Then $\beta$ restricts to an isomorphism
between  $U_{\bf r'} \setminus \beta^{-1}(P_0)$ and $U^0_{\bf r}(P_0) \setminus \{P_0\}$.
 Note that $\beta^{-1}(P_0)$ consists of the polygons of $U_{\bf r'}$ that
 degenerate at the edges $e_J$.

Define
$$g: U_{\bf r}(J^c) \rightarrow U_{\bf r'}$$ as 
$$g: (P, P_1) \rightarrow \beta^{-1}(P), \;\hbox{if}\; P \ne P_0;$$
$$g: (P_0, P_1) \rightarrow P', \;\hbox{if}\; P = P_0$$
where $P'$ is obtained from $P_1 \in M_{{\bf r}_{J^c}}$ by breaking the longest edge
(whose length is $\sum_{J^c} r_j - \epsilon_{J^c}$) into the consecutive edges
according to the lengths $r'_{j_1}, \ldots, r'_{j_s}$ (which add up to
$\sum_{J^c} r_j - \epsilon_{J^c}$, by assumption).
The inverse of $g$ can be checked to be
$$g^{-1}: P'  \rightarrow (\beta (P'), P_1),  \;\hbox{if}\; \beta (P') \ne P_0$$
where $P_1 \in M_{{\bf r}_{J^c}}$ is uniquely determined by $\beta (P')$.
$$g^{-1}: P'  \rightarrow (\beta (P'), P_1),  \;\hbox{if}\; \beta (P) = P_0$$
where $P_1$ is obtained from $P'$ (which degenerates at $e_J$) by taking the degenerating
edges as a single edge whose length is, by assumption, $\sum_J r_j' = \sum_{J^c} r_j - \epsilon_{J^c}$.
The function $g$ and $g^{-1}$ are easily seen to be continuous and holomorphic
on some dense open subsets. Thus by
the Removable Singularity Theorem (Lemma \ref{extendingmap}), $g$ and $g^{-1}$ are both holomorphic.

By the same proof as in Lemma \ref{nolinegons}, $U_{\bf r'}$ contains no line gons
(it also more or less obviously follows from the fact
that $U_{\bf r'} = \beta^{-1}(U^0_{\bf r}(P_0))$
and Lemma \ref{nolinegons}).
This shows that $U_{\bf r}(J^c) \cong U_{\bf r'}$ is smooth. Consider the projection
$$\widehat{U}_{\bf r}^1({P_0}) \rightarrow U_{\bf r}(J^c).$$
It is easy to see that all fibers of the projection are projective spaces (mostly
${\Bbb P}^0$) and the projection is bimeromorphic.
Thus the uniqueness of blowup implies that $\widehat{U}_{\bf r}^1({P_0}) \rightarrow U_{\bf r}(J^c)$
is a blowup of $U_{\bf r}(J^c)$ along smooth centers. Thus
$\widehat{U}_{\bf r}^1({P_0})$ is smooth, and 
this completes the proof.
\endproof

\begin{rem} Likewise, one obtains $U_{\bf r}(J)$ by forgetting the last factor of
$$\widehat{U}_{\bf r}^1({P_0}) \subset
U^0_{\bf r}({P_0}) \times M_{{\bf r}_J} \times M_{{\bf r}_{J^c}}.$$
By the same arguments as in the above proof,   $U_{\bf r}(J)$ can be identified with
a smooth open subset of $M_{{\bf r}''}$ for some other ${\bf r}''$. It follows then
that $U_{\bf r}(J^c) \rightarrow  U^0_{\bf r}(P_0)$
and $U_{\bf r}(J) \rightarrow  U^0_{\bf r}(P_0)$ are also resolution of
singularities, although neither is canonical due to a choice between $J$ and $J^c$.
As a consequence of the above,
 we have shown that $\widehat{U}_{\bf r}^1({P_0})$ is a common blowup of
$U_{\bf r}(J^c)$ and $U_{\bf r}(J)$.
\end{rem}

If we set $Z_J$ (compare with $Z_{[n]}$)
 to be the moduli space of $(|J|+1)$-gons with first $|J|$ sides having fixed
length $r_{j_1}, \ldots, r_{j_s}$ but the last side being free, then all the above 
discussions about $Z_{[n]}$ apply to $Z_J$ in an obvious way.

\begin{say}
\label{correpondencevariety}
This allows us to define the following correspondence (variety).
For any stable polygon
$${\bf P}=(P_0, P_1, \cdots, P_m) \in {M}_{{\bf r}} \times {M}_{{\bf r}_{J_1}} \times \cdots 
   \times {M}_{{\bf r}_{J_m}},$$
define
$$U_{\bf r}({\bf P}) = \left\{\begin{array}{ccccc}
(P_0', P_1', \cdots, P_m') \in {M}_{{\bf r}} \times {M}_{{\bf r}_{J_1}} \times \cdots 
   \times {M}_{{\bf r}_{J_m}} : \\
 \iota(P'_s, P'_t)\; \hbox{whenever}\; J_t \subset J_s
\end{array} \right\}.$$
If ${\bf P}=P_0$, we simply set $U_{\bf r}({\bf P}) = M^0_{\bf r}$.
Note that the components of a stable polygon are labeled but not ordered. Consequently,
we need to point out that the incident correspondence
$U_{\bf r}({\bf P})$ is, up to isomorphisms (induced by permutations among the
components), uniquely determined by the stable polygon ${\bf P}$.
\end{say}

\begin{lem}
Let $(P_0', P_1', \cdots, P_m') \in U_{\bf r}({\bf P})$. For every $J_t \subset J_s$,
$P_t'$ is uniquely determined by $P'_s$ unless $P_t'$ is a bubble of $P'_s$.
\end{lem}

\proof
This follows from the definition of incidence relation.
\endproof

It then follows 

\begin{cor}
There is a canonical injection 
$$\gamma_{\bf r}({\bf P}): U_{\bf r}({\bf P}) \rightarrow {\bfmit M}_{{\bf r}, \varepsilon}.$$
\end{cor}

\begin{lem}
\label{keylemma} Given any stable polygon ${\bf P}$,
$U_{\bf r}({\bf P})$ is a smooth complex analytic subvariety 
of ${M}_{{\bf r}} \times {M}_{{\bf r}_{J_1}} \times \cdots 
   \times {M}_{{\bf r}_{J_m}}$.
\end{lem}

\proof In this proof, the permutation scheme will be applied whenever the length of a diagonal
of a polygon is used.
Given any stable polygon
$${\bf P} = (P_0,  P_1, \ldots, P_m) \in  
M_{{\bf r}} \times {M}_{{\bf r}_{J_1}} \times \cdots 
   \times {M}_{{\bf r}_{J_m}},$$
we shall describe $U_{\bf r}({\bf P})$ inductively.
If $P_0$ is generic, $U_{\bf r}({\bf P}) =  M^0_{\bf r}$ is obviously complex analytic
and smooth. So we assume that
$m \ge 1$. If $J = J_k$ for some $1 \le k \le m$, we will call $P_J' := P_k'$
to be the $J$-component (or, $k$-component) of ${\bf P}'= (P_0',  P_1', \ldots, P_m')$.

First, let $U^0_{\bf r}({\bf P})$ 
be the subset of 
$M_{\bf r}$ 
consisting of points 
$P_0'$ 
with the following property:
For any $J$ such that $P_0$ degenerates at $J$,
we require that the $J$-diagonal of $P'$ satisfies
$$ \sum_J r_j - 2 \min_J \{r_j\} < |d_J| \le \sum_J r_j.$$
This is a complex analytic open neighborhood of $P_0$ in $M_{\bf r}$.

Define  an incident correspondence $\widehat{U}_{\bf r}^1 ({\bf P}) \subset 
U^0_{\bf r}({\bf P}) \times \Pi_J M_{{\bf r}_J}$
as follows. Here, the product
is over all $J$'s such that
$P_0$ degenerates at the edges $e_J$. Such $J$'s will be referred as  the subsets of the first kind. 
We then require that for any point ${\bf P}' \in 
\widehat{U}_{\bf r}^1 ({\bf P})$
we have the incidence relation $i (P_0', P_J')$ where $P_0'$ is the (main) 0-component of
${\bf P}'$ and $P_J'$ is the  $J$-component of  ${\bf P}'$.
Since the incidence relation is holomorphic,
$\widehat{U}_{\bf r}^1 ({\bf P})$ is holomorphic.

Then we take $U_{\bf r}^1({\bf P})$ to be the subset of 
$\widehat{U}_{\bf r}^1({\bf P})$ 
consisting of points with the following
additional requirements. Given any point ${\bf P}'$ in 
$\widehat{U}_{\bf r}^1({\bf P})$, 
it belongs to $U_{\bf r}^1({\bf P})$ if for any
$J$-component $P_J'$ of the point (where $J$ is of the first kind)
and any $I \subset J$ such that $P_J'$ degenerates at the edges
$e_I$, we demand that
$$ \sum_I r_j - 2 \min_I \{r_i\} < l(d_{J, I}') \le \sum_{J 
\setminus I} r_j$$
where $d_{J,I}'$ is the $I$-diagonal of $P_J'$.
Then $U_{\bf r}^1({\bf P})$ 
is a holomorphic open subset of $\widehat{U}_{\bf r}^1({\bf P})$.
We may refer to the above $I$'s  as the subsets  of the second kind.

Likewise, 
one can define the incident correspondence variety
$\widehat{U}_{\bf r}^2({\bf P}) \subset U_{\bf r}^1({\bf P}) \times 
\Pi_I M_{{\bf r}_I}$ where the product is over all subsets  of the second kind and
a complex analytic open subset $U_{\bf r}^2({\bf P}) \subset \widehat{U}_{\bf r}^2({\bf P})$, 
$\ldots$, and keep going until all $J_1, \ldots, J_m$ are taken into account.
We hence obtain inductively a sequence of projections
$$\widehat{U}_{\bf r}^1({\bf P}) \rightarrow U^0_{\bf r}({\bf P}),$$
$$\widehat{U}_{\bf r}^2({\bf P}) \rightarrow U^1_{\bf r}({\bf P}),$$
$$\ldots $$
$$\widehat{U}_{\bf r}^{h-1}({\bf P}) \rightarrow {U}_{\bf r}^{h-2}({\bf P}),$$ 
$$\widehat{U}_{\bf r}^h({\bf P}) \rightarrow {U}_{\bf r}^{h-1}({\bf P}).$$
One checks directly that the so-inductively defined
$\widehat{U}_{\bf r}^h({\bf P})$ coincides with the
incident correspondence variety  $U_{\bf r}({\bf P})$ as defined in \ref{correpondencevariety}.

To see that $U_{\bf r}({\bf P})$ is smooth, we will analyze the
above sequence. We will have to divide the proof into two cases.

{\sl Case 1}. $P_0$ is not a line gon. In this case,
$P_0$ is a smooth point of $M_{\bf r}$ and thus
the open neighborhood $U^0_{\bf r}({\bf P})$ is smooth by Lemma \ref{nolinegons}.
Every fiber of the projection $\widehat{U}_{\bf r}^1({\bf P}) \rightarrow U^0_{\bf r}({\bf P})$
 is easily seen to be projective spaces (mostly ${\Bbb P}^0$).
In addition, the projection is bimeromorphic. Thus the uniqueness
of blowup implies that the above projection is a blowup of $U^0_{\bf r}({\bf P})$
 along  smooth centers. Indeed, the centers are the loci of
$P_0' \in U^0_{\bf r}(P_0)$ that are degenerate at some edges $e_J$. Each component
of the center can be identified, by taking the degenerate edges $e_J$ as a single edge,
 with a smooth open subset of $M_{\bf x}$ (which must be a projective space because such
${\bf x}$ must be in a favorable chamber).
Thus $\widehat{U}_{\bf r}^1({\bf P})$ is smooth.
Likewise, step by step, one can show that $\widehat{U}_{\bf r}^2({\bf P})$, $\ldots$, and 
 $\widehat{U}_{\bf r}^h({\bf P})$ are all smooth. In particular,
$U_{\bf r}({\bf P})$ is smooth.

{\sl Case 2}. $P_0$ is a line gon. In this case, 
$P_0$ is a singular point of $M_{\bf r}$ and thus
the open neighborhood $U^0_{\bf r}({\bf P})$ is singular.
By Lemma \ref{canonicalresolution}, 
noting that $\widehat{U}_{\bf r}^1({\bf P}) = \widehat{U}_{\bf r}^1(P_0)$,
the projection
$$\widehat{U}_{\bf r}^1({\bf P})
\rightarrow U^0_{\bf r}({\bf P})$$ is a (canonical) resolution
of singularities of $U^0_{\bf r}({\bf P})$. 

Now, follow the proof in {\sl Case 1}, 
build up from $\widehat{U}_{\bf r}^1({\bf P})$, step by step, we conclude that
$\widehat{U}_{\bf r}^h({\bf P})$ is smooth as desired.
\endproof

As an immediate consequence of Lemma \ref{keylemma}, we have

\begin{thm}
$ {\bfmit M}_{{\bf r}, \varepsilon}$ carries a natural smooth, compact complex analytic structure induced by 
the injections $\gamma_{\bf r}({\bf P}): U_{\bf r}({\bf P}) \rightarrow {\bfmit M}_{{\bf r}, \varepsilon}$
for all stable polygons ${\bf P}$.
\end{thm}

\proof
The complex structures on $\gamma_{\bf r}({\bf P})(U_{\bf r}({\bf P}))$ 
(induced by those on $U_{\bf r}({\bf P})$) for various stable polygons ${\bf P}$
 obviously agree with each other over the overlaps. The theorem then follows.
\endproof




\section{${\bfmit M}_{{\bf r}, \varepsilon}$ as iterated blowups of $M_{\bf r}$}

Globally, by forgetting all the bubbles of a stable polygon, or locally,
by the projection from $U_{\bf r}({\bf P})$ (see 5.11) to (the main factor) $M_{\bf r}$, we obtain 

\begin{cor} There is a canonical complex analytic map
$$\pi_{{\bf r},\varepsilon} : {\bfmit M}_{{\bf r}, \varepsilon} \rightarrow M_{\bf r}$$
which restricts to the identity on $M^0_{\bf r}$.
\end{cor}

\proof
The existence of the projection $\pi_{\bf r}$ as a set-theoretic map is obvious.
That $\pi_{\bf r}$ is complex analytic follows from that locally $\pi_{\bf r}$ is equivalent
to the projection $U_{\bf r}({\bf P}) \rightarrow M_{\bf r}$.
\endproof

\begin{say}
When $M_{\bf r}$ is smooth, that is, when ${\bf r}$ is away from walls,
the map
$$\pi_{{\bf r}, {\varepsilon}}: {\bfmit M}_{{\bf r}, {\varepsilon}} \rightarrow M_{\bf r}$$
ought to be an iterated blowup of $M_{\bf r}$ along some smooth subvarieties.
The details go as follows.
\end{say}

\begin{say}
Let $F[n]$ be the set of all partitions of $[n]=\{1, \ldots, n\}$,
partially ordered by reverse refinement. That is, an element of 
 $F[n]$ is of the form $\alpha=I_1 \coprod \ldots \coprod I_k =[n]$.
Let $\beta= J_1 \coprod \ldots \coprod J_m =[n]$ be another partition
of $[n]$. We say $\alpha \le \beta$ if for any $1 \le s \le m$,
$J_s \subset I_t$ for some $1 \le t \le k$. The maximal element in
$F[n]$ is $\{1\}  \coprod \ldots \coprod \{n\}$. The smallest element is
$[n]$ itself (but we will never use this smallest element).

Given any $\alpha$ in $F[n]$, define $Y_\alpha$ to be the set of all
polygons $P$ in $M_{\bf r}$ such that the edges $e_{I_s}$ ($1 \le s \le k$)
are parallel. This is a closed subvariety of
$M_{\bf r}$ which is isomorphic to $M_{{\bf r}_\alpha}$ 
by an obvious natural map\footnote{To see this, one may need to use 
 permutations among edges. Consult Figure 7.}
 $$f_\alpha: Y_\alpha \rightarrow  M_{{\bf r}_\alpha},$$
where
$${\bf r}_\alpha =(\sum_{s \in I_1} |e_s|, \ldots, \sum_{s \in I_k} |e_s|).$$
Figure 7 illustrates the map $f_\alpha$ as the composition of a permutation
followed by an identification.

\begin{center}
\vskip   2.1in
\special{bmp:poly7.bmp x=5in y=2.1in}

Figure 7
\end{center}

We must point out that the strata $Y_\alpha$ are empty for many $\alpha
\in F[n]$ (e.g.,  when ${\bf r}_\alpha$ lies outside of the cone $C({\Bbb D}^k_2)$).

Obviously, all $Y_\alpha$ are smooth if $M_{\bf r}$ is (having no line gons).
In general, it always contains a smooth dense open subset
$$Y^0_\alpha = f_\alpha^{-1} (M_{{\bf r}_\alpha}^0).$$
That is, $Y^0_\alpha$ is the set of all
polygons $P$ in $M_{\bf r}$ such that the edges $e_{I_s}$ ($1 \le s \le k$)
are all parallel to each other  but not to any other edges.

Then we obtain a canonical decomposition
$$ M_{\bf r} = \bigcup_{\alpha \in F[n]} Y^0_\alpha.$$
This is a stratification of $ M_{\bf r}$ by smooth strata. 
\end{say}

One checks readily that

\begin{prop}
$Y_\alpha \subset Y_\beta$ if and only if $\alpha \le \beta$.
\end{prop}

In fact, it can be shown that an intersection of the closed strata $Y_\alpha$'s is again
a closed stratum unless it is empty.

Now we come to our main theorems in this section.

\begin{thm}
\label{blowups}
 Assume that $M_{\bf r}$ is smooth. Then
$$\pi_{{\bf r},\varepsilon}: {\bfmit M}_{{\bf r}, {\varepsilon}}  \rightarrow M_{\bf r}$$
is the iterated blowup of $M_{\bf r}$ along 
(the proper transforms of) all the smooth closed strata
$Y_\alpha$ in the order dictated by the partial order 
($Y_\alpha \subset Y_\beta$ if and only if $\alpha \le \beta$),
starting from the smallest ones.
\end{thm}

\proof
This follows from the details of the proof presented for Lemma 5.14
and the fact that the blowup construction is local and unique.
\endproof

\begin{rem} 
\label{symblowup}
For suitable $\varepsilon$,
${\bfmit M}_{{\bf r}, {\varepsilon}}$ comes equipped with 
a symplectic form $\Omega_{{\bf r}, {\varepsilon}}$ whose cohomology class is uniquely defined
so that $$\pi_{{\bf r},\varepsilon}: 
({\bfmit M}_{{\bf r}, {\varepsilon}},  \Omega_{{\bf r}, {\varepsilon}}) \rightarrow 
(M_{\bf r}, \omega_{\bf r})$$ can be interpreted as a symplectic blowup (see Lemma 6.44 and pp $230-231$ of
\cite{MS}). The ambiguity of the so-called symplectic blowups provides another angle to see
how natural it is the choices of ${\varepsilon}$ in the definition of stable polygons.
\end{rem}

The singular case is slightly complicated.
Assume that $M_{\bf r}$ is singular. Then $M_{\bf r}$ has isolated singularities
defined by line gons. Each of these line gons forms one of the smallest strata
in $ M_{\bf r} = \bigcup_{\alpha \in F[n]} Y^0_\alpha.$
Recall by Lemma 5.9, these singularities admit canonical resolutions.

\begin{thm}
\label{resblowups}
Assume that $M_{\bf r}$ is singular. Then 
$$\pi_{{\bf r},\varepsilon} : 
{\bfmit M}_{{\bf r}, {\varepsilon}}  \rightarrow M_{\bf r}$$
is the composite of the canonical resolutions as described in Lemma 5.9
followed by the iterated blowups of the resulting resolution
 along (the proper transforms of)
all other smooth closed strata
$Y_\alpha$ in the order dictated by the partial order 
($Y_\alpha \subset Y_\beta$ if and only if $\alpha \le \beta$),
starting from the smallest ones.
\end{thm} 

\proof
After applying  Lemma 5.9, all other arguments remain the same as for the
proof of Theorem \ref{blowups}
\endproof

\begin{rem}
When ${\bf r}$ is on a wall, $M_{\bf r}$ carries a singular K\"ahler form.
In this case, ${\bfmit M}_{{\bf r}, \varepsilon} \rightarrow M_{\bf r}$ ought to be interpreted as
a K\"ahler morphism in a suitable sense.
\end{rem}

\begin{say}
Among all the iterated blowups $$\pi_{{\bf r}, {\varepsilon}}:
({\bfmit M}_{{\bf r}, {\varepsilon}}, \Omega_{{\bf r}, {\varepsilon}})
\rightarrow (M_{\bf r}, \omega_{\bf r}),$$
two special cases worth mentioning. 
For one kind of special choices of 
${\bf r}$, $M_{\bf r}$ is isomorphic to $({\Bbb P}^1)^{n-3}$
(see \ref{rem:p1n3}). In this case,
our presentation of the blowup
${\bfmit M}_{{\bf r}, {\varepsilon}} \rightarrow M_{\bf r}$,  putting aside symplectic
structures and after showing in the next section the isomorphism between
${\bfmit M}_{{\bf r}, {\varepsilon}}$ and $\overline{\cM}_{0,n}$, specializes to
the blowup representation of $\overline{\cM}_{0,n}$,
$\overline{\cM}_{0,n}  \rightarrow ({\Bbb P}^1)^{n-3}$,
as utilized by Keel in his study on the Chow ring of
$\overline{\cM}_{0,n}$ (\cite{Keel90}). For some other  special choices of 
${\bf r}$, $M_{\bf r}$ is isomorphic to ${\Bbb P}^{n-3}$ (see \ref{propo:favorchamber}).
In this case,
our  blowup ${\bfmit M}_{{\bf r}, {\varepsilon}} \rightarrow M_{\bf r}$,
again forgetting symplectic structures,
amounts to the blowup representation of $\overline{\cM}_{0,n}$,
$\overline{\cM}_{0,n}  \rightarrow {\Bbb P}^{n-3}$,
as studied by Kapranov in \cite{Kapranov95}. We give below some details of the latter.
\end{say}

\begin{exmp}
Choose ${\bf r} \in \Delta_i$ for some fixed $1\le i \le n$ 
(see \S 3 for the characterization of $\Delta_i$). Note that a polygon $P \in 
M_{\bf r} \cong {\Bbb P}^{n-3}$
will never degenerate at the edge $e_i$. This fact together
with the identification $f_\alpha: Y_\alpha \rightarrow  M_{{\bf r}_\alpha}$ implies
that every stratum $Y_\alpha$ is isomorphic to ${\Bbb P}^{\dim Y_\alpha}$.
Point strata correspond to polygons $P$ whose edges have exactly three different directions.
There are $(n-1)$ such strata. They are points in ${\Bbb P}^{n-3}$ in general position.
Any other stratum,  obviously containing a subset  of these points, is the projective subspace
spanned by the points in the subset. Applying Theorem \ref{blowups}
(together with the isomorphism ${\bfmit M}_{{\bf r}, {\varepsilon}}\cong
\overline{\cM}_{0,n}$ to be proved in \S 7),
we recover Kapranov's blowup representation, $\overline{\cM}_{0,n}  \rightarrow {\Bbb P}^{n-3}$.
\end{exmp}

\section{${\bfmit M}_{{\bf r},\varepsilon}$ and $\overline{\cM}_{0,n}$}

\begin{say}
Recall that a $n$-pointed connected complex algebraic curve of genus $0$ is stable if
\begin{enumerate}
\item the $n$-marked points are smooth points;
\item  every singular point is an ordinary double point;
\item for each irreducible connected component, the number of  marked points plus
      the number of singular points on the component is at least 3.
\end{enumerate}
The set $\overline{\cM}_{0,n}$
of equivalence classes of all $n$-pointed stable algebraic curve of genus $0$ carries a natural
structure of a smooth projective variety. 
Let $\cM_{0,n}$ be the moduli space of 
$n$-pointed smooth curves of genus zero.
Then $\overline{\cM}_{0,n} \setminus \cM_{0,n}$ is a divisor with normal crossings.
Given any $n$-pointed stable curve $X$,
we can associate to it a graph $\graph(X)$: the vertices of $\graph(X)$ correspond
to the irreducible components of $X$, and two  vertices are joined by an edge if
their corresponding components share a common singular point. Then that the curve
$X$ is of genus $0$ is equivalent to  the graph $\graph(X)$ being a tree.
\end{say}

\begin{say}
Likewise, we can also attach a graph $\graph({\bf P})$ to any given stable
$n$-gon ${\bf P}=(P_0, P_1 \ldots, P_m)$: the vertices of $\graph(X)$ correspond
to the polygons $\{P_0, P_1 \ldots, P_m\}$, and two vertices are joined by an edge if
their corresponding polygons satisfy the bubble relation (Definition \ref{bubblerelation}).
One checks easily that $\graph({\bf P})$ is a tree.
\end{say}

\begin{thm}
\label{isomorphismtoM_{0,n}-bar}
Let ${\bf r}$ be a point in the interior of $C({\Bbb D}^n_2)$ and ${\varepsilon}$ be chosen
as before. Then,
${\bfmit M}_{{\bf r},{\varepsilon}} $ 
and $\overline{\cM}_{0,n}$ are biholomorphic.
Consequently, the complex structure on ${\bfmit M}_{{\bf r},{\varepsilon}}$ is independent\footnote{
The structure of K\"ahler space, however, depends on choices.} of
${\bf r}$ and the choice of ${\varepsilon}$. Moreover, 
  ${\bfmit M}_{{\bf r},{\varepsilon}}\setminus M^0_{\bf r} $ is a divisor with
normal crossings.
\end{thm}

\proof 
${\bfmit M}_{{\bf r}, \varepsilon}$ is easily seen to be bimeromorphic to
$\overline{\cM}_{0,n}$. That is, we have a biholomorphism
$$\gamma^0: M^0_{\bf r} (\subset {\bfmit M}_{{\bf r},{\varepsilon}})
 \rightarrow \cM_{0,n}  (\subset \overline{\cM}_{0,n}).$$
This map can be extended continuously as follows.
For any stable polygon 
$${\bf P}=(P_0, P_1, \ldots, P_m) \in {M}_{{\bf r}} \times {M}_{{\bf r}_{J_1}} \times \cdots 
   \times {M}_{{\bf r}_{J_m}},$$
we obtain a collection 
$X=(X_1, \ldots, X_m)$ of pointed curves of genus zero
via the identifications
$$M_{\bf r} \cong ({\Bbb P}^1)^n_{ss}({\bf r}_C)/\!/\PGL (2).$$
Now regarding the coinciding points in each $X_i$ as a single point,
and if $(P_s, P_t)$ is a bubble pair, we joint $X_s$ and
$X_t$ at the coinciding points (regarded as single point)  of $X_s$ 
(which corresponds to the same direction pointing edges of the polygon $P_s$)
 and the point of $X_t$ that corresponds to the longest edge of $P_t$.
This way, we obtain a reducible algebraic curve $X_1 \cup \ldots \cup X_m$
which has $n$-labeled (smooth) points. It is of genus zero because
its associated graph coincides with the graph of the stable polygon ${\bf P}$
which is a tree by the construction. One checks that each component $X_i$
has at least three distinguished (marked plus singular) points. Hence
$X_1 \cup \ldots \cup X_m$ is a stable $n$-pointed curve of genus zero.
The so-induced map 
$$\gamma: {\bfmit M}_{{\bf r}, \varepsilon} \rightarrow \overline{\cM}_{0,n}$$
$$(P_1, \ldots,P_m) \to X_1 \cup \ldots \cup X_m$$
is (more or less) obviously continuous and injective.
By Lemma \ref{extendingmap},
$\gamma$ is holomorphic.
Since $\gamma$ is bimeromorphic, it is also surjective.
This implies that $\gamma$ is biholomorphic.
The rest follows from some well-known properties of $\overline{\cM}_{0,n}$.
\endproof

\begin{rem}
When ${\bf r}$ is away from walls, being K\"ahler and Moishezon, 
${\bfmit M}_{{\bf r}, \varepsilon}$ is projective. Since
${\bfmit M}_{{\bf r}, \varepsilon}$ and $\overline{\cM}_{0,n}$ are biholomorphic,
the GAGA theorems imply that
${\bfmit M}_{{\bf r}, \varepsilon}$ and $\overline{\cM}_{0,n}$ are isomorphic
as projective varieties as well.
\end{rem}


\section{The K\"ahler cone of $M_{\bf r}$}

\begin{say}
To pave a way to determining the K\"ahler cone of $\overline{\cM}_{0,n}$,
we first study the K\"ahler cone of $M_{\bf r}$. We shall mainly  focus
our attention on  the most important
special cases that ${\bf r}$ lies in a chamber around the center of $C({\Bbb D}^n_2)$.
The others, though less significant for the K\"ahler cone of $\overline{\cM}_{0,n}$,
can be treated similarly. There will be some differences between the case when
$n$ is odd and the case when $n$ is even.
\end{say}

\begin{say} {\sl The odd case}. This is the nicer case. There is a unique chamber
$C_0$ that contains the half ray ${\Bbb R}_+ \cdot (1, \ldots, 1)$ in $C({\Bbb D}^n_2)$.
In terms of inequality, $C_0$
is defined by
$$ 2 \sum_{j \in J} r_j < \sum_{1 \le j \le n} r_j\; \hbox{if} \; |J| < {n \over 2};$$
$$ 2 \sum_{j \in J} r_j > \sum_{1 \le j \le n} r_j\; \hbox{if} \; |J| > {n \over 2}.$$
Equivalently, $C_0 \cap{\Bbb D}^n_2 $ is defined by
$$  \sum_{j \in J} r_j < 1\; \hbox{if} \; |J| < {n \over 2};$$
$$ \sum_{j \in J} r_j > 1 \; \hbox{if} \; |J| > {n \over 2}.$$

Recall that
all  $M_{\bf r}$ underly a common
projective variety $M_{C_0}$
when ${\bf r}$ lies in the interior of $C_0$.
It is well known that the second Betti number
or the Picard number of this projective variety is equal to $n$.
Many people have made independent calculations. The following provides one of them.


\begin{thm} 
\label{odd}
$M_{C_0}$ can be obtained 
from ${\Bbb P}^{n-3}$ by the following sequence of ``blowups'' and ``blowdowns'':
we start with blowing up ${n-1\choose 1}$ points (in general position) of ${\Bbb P}^{n-3}$
to ${n-1\choose 1}$ ${\Bbb P}^{n-4}$; then blowing
down ${n-1\choose 2}$ ${\Bbb P}^1$ to ${n-1\choose 2}$ points followed by blowing up 
${n-1\choose 2}$ points to ${n-1\choose 2}$ ${\Bbb P}^{n-5}$;
then blowing
down ${n-1\choose 3}$ ${\Bbb P}^2$ to ${n-1\choose 3}$ points followed by blowing up 
${n-1\choose 3}$ points to ${n-1\choose 3}$ ${\Bbb P}^{n-6}$; ...;
finally, blowing down ${n-1\choose \frac{n-3}{2}}$ ${\Bbb P}^{\frac{n-5}{2}}$ to 
${n-1\choose \frac{n-3}{2}}$ points
 followed by blowing up 
${n-1\choose \frac{n-3}{2}}$ points to ${n-1\choose \frac{n-3}{2}}$
${\Bbb P}^{\frac{n-3}{2}}$. 
\end{thm}

\proof (Outline.) Take a general point in a favorable chamber of ${\Bbb D}^n_2$.
 Consider the line segment joining this point and the
  barycentre of ${\Bbb D}^n_2$. By Theorem 2.2 of \cite{Hu91},
   one get a sequence of projective morphisms between 
  ${\Bbb P}^{n-3}$ and $M_{C_0}$
such that the dimensions of special fibers change according to the rule
$$d+e=n-4,$$
starting with $0+(n-4)$, then $1+(n-5)$, and so on ..., till we arrive at
the center of ${\Bbb D}^n_2$ and end up with ${n-5 \over 2} + {n-3 \over 2}$.
Taking into account of
the number of walls that we have crossed, we obtain the assertion as in the theorem.
\endproof

The theorem immediately implies that
{\footnotesize{
$${\bf P}_t(M_{C_0}) ={\bf P}_t({\Bbb P}^{n-3}) +  {n-1\choose 1} [{\bf P}_t({\Bbb P}^{n-4})-1]
+ {n-1\choose 2} [{\bf P}_t({\Bbb P}^{n-5}) -{\bf P}_t({\Bbb P}^1)] $$
$$+ \cdots +
 {n-1\choose\frac{n-3}{2}} [ {\bf P}_t({\Bbb P}^{\frac{n-3}{2}}) -{\bf P}_t({\Bbb P}^{\frac{n-5}{2}})].$$ }}
That is,
{\footnotesize{
$${\bf P}_t(M_{C_0}) = {{t^{2(n-2)} -1} \over {t^2 -1}} + {n-1\choose 1} t^2 
{{t^{2(n-4)} -1} \over {t^2 -1}} + {n-1\choose 2} t^4 {{t^{2(n-6)} -1} \over {t^2 -1}}$$
$$+ \cdots + {n-1\choose\frac{n-3}{2}} t^{n-3}.$$ }}
Here
${\bf P}_t(X)$ denotes the Poincar\'e polynomial of $X$.
These  Betti numbers were probably first computed by Kirwan.

Now we turn to the K\"ahler cone of $M_{C_0}$.
Define a map
$$\theta_0: C_0 \rightarrow K(M_{C_0})$$ by
$$ r \to [\omega_{\bf r}],$$
where $[\omega_{\bf r}]$ is the cohomology class of $\omega_{\bf r}$.
This is initially defined in the interior of $C_0$ but can be extended easily
to the boundary. 
\end{say}

\begin{thm}
\label{n=odd}  Assume $n \ge 5$.
The  K\"ahler cone $K(M_{C_0})$ of $M_{C_0}$ can be naturally identified with the cone
$C_0$ by the linear isomorphism $\theta:  C_0 \rightarrow K(M_{C_0})$.
\end{thm}

\proof  By Duistermaat-Heckman's theorem\footnote{
Strictly speaking, to apply Duistermaat-Heckman's theorem in this case,
one should use the Grassmannian construction of
$M_{\bf r}$ as mentioned in \ref{grasscons}.
}, $\theta$ is linear.
Since $\dim C_0 = \dim K(M_{C_0})=n$,
to prove that $\theta$ is a linear 
isomorphism, it suffices to argue that $\theta_0$ is surjective
onto the open K\"ahler cone. Because the K\"ahler cone is generated
by the cone of ample divisors (in this case), it is enough to consider integral
points. Let $L$ be any ample line bundle over $M_{C_0} \cong ({\Bbb P}^1)^n_{ss}({\bf r}_{C_0})
/\!/ GL(2, {\Bbb C})$ where ${\bf r}_{C_0}$ is any fixed integral point in
the interior of $C_0$. Let $\pi: ({\Bbb P}^1)^n_{ss}({\bf r}_{C_0})
\rightarrow ({\Bbb P}^1)^n_{ss}({\bf r}_{C_0})
/\!/ GL(2, {\Bbb C})$ be the quotient map. Then
$\pi^* L$ extends canonically to an ample line bundle $\cL$ over $({\Bbb P}^1)^n$
because the complement $({\Bbb P}^1)^n \setminus ({\Bbb P}^1)^n_{ss}({\bf r}_{C_0})$
has codimension greater than 1. Now the surjectivity follows from
the fact that $({\Bbb P}^1)^n_{ss}({\bf r}_{C_0}) = ({\Bbb P}^1)^n_{ss} (\cL)$
and $\cL$ descends to $L$.
It remains to show that the image $\theta (\partial C_0)$ lies on the boundary
of $K(M_{C_0})$. This basically follows from the fact that when
${\bf r}$ approaches the boundary of $C_0$, one gets a non-trivial (symplectic)
blow down map, and in particular, for any curve $S$ that is contracted to a point the pairing
$\langle [\omega_{\bf r}], [S] \rangle$ approaches zero as 
${\bf r}$ approaches the boundary of $C_0$. Here $[S]$ is the homology class of $S$.
\endproof

\begin{rem} A few remarks are in order. First,
$C_0$ is the unique chamber that is invariant under permutations of coordinates.
Second, $M_{C_0}$ is the only quotient that admits an induced
       action of the permutation group $\Sigma_n$ (which acts on $({\Bbb P}^1)^n$ by
       permuting the coordinates). This is so because $({\Bbb P}^1)^n_{ss}({\bf r}_{C_0})$
is the only semi-stable set that is invariant under $\Sigma_n$.
\end{rem}

The even cases are somewhat complicated due to the fact that the half ray 
$E={\Bbb R}_+ \cdot (1, \ldots, 1)$ lies in the intersection of a number of
chambers. 

First we record a result on $M_E$.

\begin{thm} 
$M_E$ can be obtained 
from ${\Bbb P}^{n-3}$ by the following sequence of ``blowups'' and ``blowdowns'':
 we start with blowing up ${n-1 \choose 1}$ points
to ${n-1\choose 1}$ ${\Bbb P}^{n-4}$; then blowing
down ${n-1\choose 2}$ ${\Bbb P}^1$ to ${n-1\choose 2}$ points followed by blowing up 
${n-1\choose 2}$ points to ${n-1\choose 2}$ ${\Bbb P}^{n-5}$;
then blowing
down ${n-1\choose 3}$ ${\Bbb P}^2$ to ${n-1\choose 3}$ points followed by blowing up 
${n-1\choose 3}$ points to ${n-1\choose 3}$ ${\Bbb P}^{n-6}$; ...;
then blowing
down ${n-1\choose \frac{n-4}{2}}$ ${\Bbb P}^{\frac{n-6}{2}}$ to 
${n-1\choose\frac{n-4}{2}}$ points followed by blowing up 
${n-1\choose\frac{n-4}{2}}$ points to ${n-1\choose\frac{n-4}{2}}$
${\Bbb P}^{\frac{n-2}{2}}$; finally, blowing down 
${n-1\choose\frac{n-2}{2}}$ ${\Bbb P}^{\frac{n-4}{2}}$ to 
${n-1\choose\frac{n-2}{2}}$ points. 
\end{thm}

\proof
The proof is 
similar to that of Theorem \ref{odd}.
\endproof

It follows that
{\footnotesize{
$${\bf I}\!{\bf P}_t(M_E) ={\bf P}_t({\Bbb P}^{n-3}) +  {n-1\choose 1} [{\bf P}_t
({\Bbb P}^{n-4})-1]
+ {n-1\choose 2} [{\bf P}_t({\Bbb P}^{n-5}) -{\bf P}_t({\Bbb P}^1)] $$
$$+ \cdots +
 {n-1\choose \frac{n-4}{2}} [ {\bf P}_t({\Bbb P}^{\frac{n-2}{2}}) -{\bf P}_t({\Bbb P}^{\frac{n-6}{2}})].$$}}
That is,
{\footnotesize{
$${\bf I}\!{\bf P}_t(M_{E}) = {{t^{2(n-2)} -1} \over {t^2 -1}} + {n-1\choose 1} t^2 
{{t^{2(n-4)} -1} \over {t^2 -1}} + {n-1\choose 2} t^4 {{t^{2(n-6)} -1} \over {t^2 -1}}$$
$$+ \cdots + {n-1\choose\frac{n-4}{2}} t^{n-4} {{t^4-1} \over {t^2-1}}.$$}}
Here ${\bf I}\!{\bf P}_t(X)$ denotes the intersection  Poincar\'e polynomial of $X$ (\cite{GM}).
These intersection Betti numbers are first computed by Kirwan.

Now, take any maximal chamber $C_*$ that contains the half ray $E$
(the symbol $*$ stresses on the choice of a base point to single out a chamber).
Again, the Picard number of $M_{C_*}$ is equal to $n$, and we can define a map
$$\theta: C_* \rightarrow K(M_{C_*}),$$
in the same way as before. A similar proof as in Theorem \ref{n=odd} will give

\begin{thm}
\label{n=even}  Assume $n \ge 6$.
The  K\"ahler cone $K(M_{C_*})$ of $M_{C_*}$ can be naturally identified with the cone
$C_*$ by the linear isomorphism $\theta:  C_* \rightarrow K(M_{C_0})$.
\end{thm}


\section{The K\"ahler cone of $\overline{\cM}_{0,n}$}

\begin{say}
By the virtue of Remark \ref{symblowup}, $\overline{\cM}_{0,n}$ carries (classes of) symplectic forms
that are transported from (the classes of) $\Omega_{{\bf r}, \varepsilon}$
by the isomorphism $\gamma: {\bfmit M}_{{\bf r}, \varepsilon} \rightarrow \overline{\cM}_{0,n}$. 
A quick computation on the number of independent parameters for $\Omega_{{\bf r}, \varepsilon}$
leads to $n+ {{2^n-2-2n-n(n-1)}\over{2}} = 2^{n-1}- {{n^2-n+2}\over{2}}$ which coincides with the
second Betti number or the Picard number of $\overline{\cM}_{0,n}$. It is well-known that
the K\"ahler cone (or dually the Mori cone of effective curves)
 of a projective variety is,  in general,  very hard to compute. 
Our theory, as providing a large family of (classes of) K\"ahler forms 
$\Omega_{{\bf r}, \varepsilon}$, sheds a light on the shape of the K\"ahler cone of $\overline{\cM}_{0,n}$.
We shall now give some heuristic arguments and formulate a conjecture below. 
\end{say}

\begin{say}
To this end, let $C^*$ be $C_0$ when $n$ is odd or
a choice $C_*$ of the maximal chambers that contains the ray ${\Bbb R}_+ \cdot (1, \ldots, 1)$
when $n$ is even, and introduce a large cone built on 
$${\bf r} \in C^* \;\hbox{and}\; 
\varepsilon=(\epsilon_J)_{J \in \cR_{>2}({\bf r})} \in {\Bbb R}_+^{2^{n-1}-1-n-{n(n-1)\over 2}}$$ 
in the legal range as follows:
$$\cC = \{({\bf r}, \varepsilon) \in C^* \times
{\Bbb R}_+^{2^{n-1}-1-n-{n(n-1)\over 2}}
| 0 < \epsilon_J < 2 \min \{r_j\}_{j \in J}, J \in \cR_{>2} ({\bf r}) \}.$$
(One can check that this is indeed a positive convex cone.) 
Assume that we have obtained a K\"ahler form
$\Omega_{{\bf r}, \varepsilon}$ for suitable $\varepsilon$, provided that ${\bf r}$ is away from walls.
 This leads to a well-defined map $\Theta^0$
from an open subcone $\cC^0$ of $\cC$ to the  K\"ahler cone $K(\overline{\cM}_{0,n})$  of  $\overline{\cM}_{0,n}$
by taking the cohomology classes of $\Omega_{{\bf r}, \varepsilon}$.
One checks that $\Omega_{{\bf r}, \varepsilon}$ depends on $({\bf r}, \varepsilon)$ 
 continuously and $\Theta^0$
extends to a continuous map
$$\Theta: \cC^0 \rightarrow K(\overline{\cM}_{0,n}).$$
To exclude triviality, we assume that $n \ge 5$. 
\end{say}

\begin{prob} Assume that $n \ge 5$. Describe $\cC^0$ and prove that
the map $\Theta$ identifies the cone $\cC^0$ with
the  K\"ahler cone $K(\overline{\cM}_{0,n})$
of $\overline{\cM}_{0,n}$.
\end{prob}




Dually, we may also consider the Mori cone of
effective curves (see a conjecture by Fulton as formulated in
\cite{KMc}).  It should be instructive to study the two approaches all together.

We end our exposition by a digressive remark.

\begin{rem}
\label{amplecones}
By 3.3.21 of \cite{DH}, $C({\Bbb D}^n_2)$ is the $G$-ample cone for
both the $\PGL (2)$-action on $({\Bbb P}^1)^n$ and the maximal torus action
on the Grassmannian $G(2, {\Bbb C}^n)$. It is known that the Chow quotients of
these two actions can be identified with $\overline{\cM}_{0,n}$. Thus the above
conjecture would establish an interesting connection between
 the $G$-ample cone of a projective $G$-variety and
the ample cone of its Chow quotient. This and moreover
the case for a general algebraic  group action
 call for further investigation.
\end{rem}

In a forthcoming paper,
we will return to these topics.

\bibliographystyle{amsplain}
\makeatletter \renewcommand{\@biblabel}[1]{\hfill#1.}\makeatother

\vskip .3cm



\noindent
Department of Mathematics, University of Texas, Arlington, TX 76019.
hu@@math.uta.edu

\end{document}